\documentclass[]{article}

\usepackage{amsfonts,amsmath,amssymb}
\usepackage{amsthm}
\usepackage[cp866]{inputenc}
\usepackage[english]{babel}
\usepackage[abs]{overpic}

\allowdisplaybreaks \theoremstyle{definition}
\newtheorem{Definition}{Definition}

\theoremstyle{plain}
\newtheorem{Theorem}{Theorem}
\newtheorem{Corollary}{Corollary}
\newtheorem{Proposition}{Proposition}

\newcommand{\bfSig}{\mathbf{\Sigma}}

\newcommand{\calN}{\mathcal{N}}

\newcommand{\bdom}{\omega^{\leq\omega}}
\newcommand{\str}{\omega^{<\omega}}
\newcommand{\cdom}{n^{\leq\omega}}
\newcommand{\conc}{\hat{\omega}^{\omega}}
\newcommand{\home}{\hat{\omega}}

\newcommand{\bo}{\mathbf{0}}

\newcommand{\bi}{\mathbf{i}}
\newcommand{\boo}{\mathbf{1}}

\newcommand{\bu}{\mathbf{u}}

\newcommand{\bff}{\mathbf{f}}

\title{Non-Collapse of the Effective Wadge Hierarchy}

\author{Victor  Selivanov
\\A.P. Ershov
Institute of
Informatics Systems SB RAS\\ 
{\tt vseliv@iis.nsk.su}
}

\begin{document}
\date{}

\maketitle

\begin{abstract}
We  study  the recently suggested effective  Wadge hierarchy in effective spaces, concentrating  on the non-collapse property. Along with hierarchies of sets, we  study  hierarchies of $k$-partitions  which are interesting on their own. In particular, we establish  sufficient conditions for the non-collapse of the effective Wadge hierarchy  and apply them to some concrete spaces.

\end{abstract}

%
%

\section{Introduction}\label{in}

Hierarchies are basic tools for calibrating objects according to their complexity, hence the non-collapse of  a natural hierarchy is fundamental for understanding the corresponding notion of complexity. A lot of papers investigate the non-collapse property in different contexts, see  e.g. \cite{s08a} for a survey of hierarchies relevant to those studied in this paper.

The Wadge hierarchy (WH), which is fundamental for descriptive set theory (DST), was  developed for the Baire space $\calN$, first for the case of sets \cite{wad84}, and recently for the $Q$-valued Borel functions on $\calN$, for any better quasiorder $Q$ \cite{km19}. A convincing extension of the $Q$-Wadge hierarchy to arbitrary topological spaces was developed in \cite{s19a}  (see also \cite{pe15,s17}).
In \cite{s19} we introduced and studied the effective Wadge hierarchy (EWH) in effective spaces as an instantiation  of the fine hierarchy (FH) \cite{s08a}. Here we concentrate on the  non-collapse property of the EWH.
As  in \cite{s19}, along with the EWH of sets we consider the EWH of $k$-partitions for $k>2$ (sets correspond to 2-partitions). 

The non-collapse of EWH is highly non-trivial already for the discrete space $\mathbb{N}$ of natural numbers. In fact, for the case of sets it follows from the results on the non-collapse of the FH of arithmetical sets introduced  in \cite{s83}; $m$-degrees of complete sets in levels of this hierarchy are among the ``natural $m$-degrees'' studied recently in  \cite{km18}. 
The case of   $k$-partitions with $2<k\leq\omega$ was also considered in \cite{s83} where many useful technical facts were obtained, but the non-collapse property was not established and even formulated because that time we did not have a convincing notion of a hierarchy of $k$-partitions (introduced only in \cite{s12,s17}). In the present paper, we prove some additional facts which, together with the results  in \cite{s83}, imply the non-collapse of EWH of $k$-partitions in $\mathbb{N}$. We  show that all levels of this hierarchy have complete $k$-partitions which are also natural in the sense of   \cite{km18}. 

We also prove the non-collapse of  EWH in  $\calN$, providing an effective version for the fundamental result in \cite{km19}; modulo this result, our proofs for $\calN$ are easy. 
Along with the spaces $\mathbb{N}$ and $\calN$, which are central in  computability theory, we discuss the non-collapse  of EWH for  other spaces which became popular in computable analysis and effective DST. The preservation property of the EWH established in \cite{s19a,s19} implies that the non-collapse property is inherited by the (effective) continuous open surjections which suggests a method for proving non-collapse. Unfortunately, this method is less general than the dual inheritance method for the Hausdorff-Kuratowski property \cite{s19a,s19}, that completely reduces this property in (computable) quasi-Polish spaces to that in the Baire space. Nevertheless, the method suggested here provides some insight which enables e.g. to show that the non-collapse property is hard to prove for the majority of spaces.

The FH (as also most of objects related to the WH) has  inherent combinatorial complexity resulting in rather technical notions and  involved proofs.  For this reason, it was not possible to make this paper completely self-contained. To make it more readable, we cite several known results and recall the most principal definitions, often with new observations and additional details. Perhaps, these efforts still do not make this paper  self-contained but, with  papers \cite{s83,s12,s19,km19} at hand, the reader would have everything to understand  the remaining technical details.  

After some preliminaries  in the next section, we recall  in Section \ref{wkpart} necessary information on the EWH. In Section \ref{ncol} we define some versions of the non-collapse property and relate them to the preservation property. While for hierarchies of sets the non-collapse property is defined in an obviousl way, for hierarchies of $k$-partitions with $k>2$ the situation is different and we meet some variants of non-collapse which are not visible in the case $k=2$. We carefully define these versions and prove relations between them, as well as their preservation properties. Hopefully, this provides  tools which could be of use for investigating the non-collapse property in other spaces.

Section \ref{app} contains main technical results of this paper which, as already mentioned, completes previous partial results from \cite{s83}. These results are purely computability theoretic and do not use topological notions at all. This provides for the readers interested in computability but not in topology the option to escape the topological part of the paper, by reading only Sections \ref{prel}, \ref{wkpart}, \ref{ncol} (restricted to the FH of $k$-partitions over the arithmetical hierarchy), and \ref{app}. 

In Section \ref{ncol1} we establish the non-collapse of the EWH of $k$-partitions in the Baire space, the domain  of finite and infinite strings, and some related spaces. Although proofs here are short (due to  the possibility to refer to some notions and proofs in \cite{km19}), the formulations are new and hopefully interesting because they provide (along with the results in Section \ref{app}) effective versions for the results in \cite{km19}.

\section{Preliminaries}\label{prel}

In this section we recall some notation, notions and facts used throughout th paper. We use standard set-theoretical notation, in particular, $Y^X$ is the set of functions from $X$ to $Y$, and $P(X)$ is the class of subsets of a set $X$.

\subsection{Effective hierarchies}\label{efhier}

Here  we briefly recall some notation and terminology about effective hierarchies in effective spaces. More details may be found e.g. in \cite{s15}. 

 All (topological) spaces  in this paper are  countably based $T_0$ (cb$_0$-spaces, for short). 
An {\em  effective cb$_0$-space} is a pair $(X,\beta)$ where $X$ is a cb$_0$-space, and $\beta:\omega\to P(X)$ is a numbering of a base in $X$ such that there is a uniformly c.e. sequence $\{A_{ij}\}$ of c.e. sets with
$\beta(i)\cap\beta(j)=\bigcup\beta(A_{ij})$ where $\beta(A_{ij})$ is the image of $A_{ij}$ under $\beta$. We simplify $(X,\beta)$ to $X$ if $\beta$ is clear from the context. The \emph{effectively open sets} in $X$ are the sets~$\bigcup\beta(W)$, for some c.e.~set~$W\subseteq\mathbb{N}$. The standard numbering $\{W_n\}$ of c.e. sets \cite{ro67} induces a numbering of the effectively open sets. The notion of effective cb$_0$-space allows to define e.g. computable and effectively open functions between such spaces \cite{wei00,s15}.

Among effective cb$_0$-space are: the discrete space $\mathbb{N}$ of natural numbers, the Euclidean spaces ${\mathbb R}^n$, the Scott domain $P\omega$ (see \cite{aj} for information about domains), the Baire space $\mathcal{N}=\mathbb{N}^\mathbb{N}$, the Baire domain $\bdom$  of  finite and infinite strings over $\omega$ with the Scott topology, the Cantor space  $2^\omega$ of binary infinite strings, the Cantor domains $\cdom$, $2\leq n<\omega$, of finite and infinite strings over $\{0,\ldots,n-1\}$ with the Scott topology; all these spaces come  with natural numberings of bases. The space $\mathbb{N}$ is trivial  topologically but very interesting for computability theory. 

Quasi-Polish spaces (introduced in \cite{br}) are important for DST and have several  characterisations.
Effectivizing one of them we obtain the following notion  identified implicitly in \cite{s15} and explicitly in \cite{br1,hs}:
 a {\em computable quasi-Polish space}  is an effective cb$_0$-space  $(X,\beta)$ such that there exists a computable effectively open surjection  from $\calN$ onto $(X,\beta)$. Most spaces of interest for computable analysis and effective DST, in particular the aforementioned ones, are  computable quasi-Polish.

Effective hierarchies of sets were studied by many authors, see e.g. \cite{mo09,bra05,bg12,hem06,s15}. We use standard $\Sigma,\Pi,\Delta$-notation, with suitable indices, for notation of levels of hierarchies of sets.
Let  $\{\Sigma^0_{1+n}(X)\}_{n<\omega}$ be the effective Borel hierarchy (we consider in this paper mostly finite levels), and  $\{\Sigma^{-1,m}_{1+n}(X)\}_n$ (with $\Sigma^{-1,1}$ usually simplified to $\Sigma^{-1}$) be the effective Hausdorff difference hierarchy over $\Sigma^0_m(X)$  in arbitrary effective cb$_0$-space $X$. Note that $\{\Sigma^0_{1+n}(\mathbb{N})\}_{n<\omega}$ coincides with the arithmetic hierarchy. We do not repeat  standard definitions but mention that the effective hierarchies come with standard numberings of all levels, so we can speak e.g. about uniform sequences of sets in a given level. We use definitions based on set operations (see e.g.  \cite{s15}); there is also an equivalent approach based on the Borel codes \cite{lo78,ke95}.  E.g., $\Sigma^0_1(X)$ is the class of effectively open sets in $X$,   $\Sigma^{-1}_2(X)$ is the class of differences of $\Sigma^0_1(X)$-sets, and $\Sigma^0_2(X)$ is the class of effective countable unions of $\Sigma^{-1}_2(X)$-sets. 

Levels of  effective hierarchies are denoted in the same manner as levels of the
corresponding classical hierarchies, using the lightface letters
$\Sigma,\Pi,\Delta$ instead of the boldface
$\mathbf{\Sigma},\mathbf{\Pi},\mathbf{\Delta}$ used for the classical
hierarchies \cite{ke95,mo09}. 

In fact, any lightface notion in this paper will have a classical boldface counterpart, as is standard in DST. In particular, recall that a function $f:X\to Y$ is {\em $\bfSig^0_2$-measurable} if
$f^{-1}(B)\in\bfSig^0_2$ for each $B\in \bfSig^0_1(Y)$. The corresponding lightface version is as follows: A function $f:X\to Y$ between effective cb$_0$-spaces is {\em $\Sigma^0_2$-measurable} if
$f^{-1}(B)\in\Sigma^0_2$ for each $B\in \Sigma^0_1(Y)$, effectively on the indices.  Every $\Sigma^0_2$-measurable function is $\bfSig^0_2$-measurable. 

The {\em Wadge reducibility} on subsets of a space $X$ is the many-one reducibility by continuous functions on $X$; it is denoted by $\leq_W$ (or, more precisely, by $\leq_W^X$). The {\em effective Wadge reducibility} on subsets of an  effective cb$_0$-space $X$ is the many-one reducibility by computable functions on $X$; it is denoted by $\leq_{eW}$ (or, more precisely, by $\leq_{eW}^X$). The  effective Wadge reducibility on subsets of $\mathbb{N}$ coincides with $m$-reducibility. The (effective) Wadge reducibility is extended to $k$-partitions $A,B$ is a straightforward way: a function $f$ on $X$ reduces $A$ to $B$, if $A=B\circ f$.

\subsection{Preorders and semilattices}\label{notation}

We use standard set-theoretical notation. In particular, $Y^X$ is the set of functions from $X$ to $Y$,
$P(X)$ is the class of subsets of a set $X$,   $\overline{C}=X\setminus C$ is the complement  $C\subseteq X$ in $X$. The domain and range of a function $f$ are denoted respectively by $dom(f)$ and $rng(f)$.
A class  $\mathcal{C}\subseteq P(X)$ has the {\em reduction property}  if for any $C_0,C_1\in\mathcal{C}$   there are disjoint $R_0,R_1\in\mathcal{C}$ such that $R_0\subseteq C_0,R_1\subseteq C_1$, and $R_0\cup R_1=C_0\cup C_1$. Note that if $\mathcal{C}\subseteq P(X)$ has the  reduction property and is closed under finite unions and intersections then  for any finite sequence $C_0,\ldots,C_n\in\mathcal{C}$   there are pairwise disjoint $R_0\ldots,R_n\in\mathcal{C}$ such that $R_i\subseteq C_i$ for every $i\leq n$, and $R_0\cup\cdots\cup R_n=C_0\cup\cdots\cup  C_n$; we call such  $(R_0\ldots,R_n)$  a {\em reduct} of $(C_0,\ldots,C_n)$.

We assume the reader to be familiar with the standard terminology and notation related to parially ordered sets (posets) and preorders.  Recall that a {\em semilattice} is a structure $(S;\sqcup)$ with binary operation $\sqcup$ such that $(x\sqcup y)\sqcup z=x\sqcup (y\sqcup z)$, $x\sqcup y= y\sqcup x$  and $x\sqcup x=x$, for all $x,y,z\in S$. By $\leq$ we denote the induced partial order on $S$: $x\leq y$ iff $x\sqcup y=y$. The operation $\sqcup$ can be recovered from $\leq$ since $x\sqcup y$ is the supremum of $x,y$ w.r.t. $\leq$. The semilattice is {\em distributive} if $x\leq y\sqcup z$ implies that $x= y'\sqcup z'$ for some $y'\leq y$ and $z'\leq z$. All semilattices considered in this paper are distributive (sometimes after adjoining a new smallest element denoted by $\bot$). A semilattice $(S;\sqcup,\leq)$ is a {\em d-semilattice} if it becomes distributive after adjoining to $S$ a new smallest element $\bot$; below we usually abbreviate $S\cup\{\bot\}$ to $S_\bot$.

An element $x$ of the semilattice $S$ is {\em join-reducible} if it  can be represented as the  supremum of two elements strictly below $x$. Element $x$  is {\em join-irreducible} if it is not join-reducible.  We denote by $I(S;\sqcup,\leq)$ the set of non-smallest join-irreducible elements of a semilattice $(S;\sqcup,\leq)$. If $S$ is distributive then $x$  is  join-irreducible iff $x\leq y\sqcup z$ implies that $x\leq y$ or $x\leq z$. By a {\em decomposition of a non-smallest element $x$} we mean a representation $x=x_0\sqcup\cdots\sqcup x_n$ where the {\em components $x_i$} are join-irreducible and pairwise incomparable (a smallest element may be considered as the supremum of the empty family of join-irreducible elements). Such a decomposition is {\em canonical} if it is unique up to a permutation of the components. Clearly, if $S$ 
 is a well founded semilattice then any non-smallest element $x\in S$ has a decomposition, and if $S$ is distributive then $x$ has a canonical decomposition.
 
To simplify notation, we often apply the terminology about posets to preorders meaning the corresponding quotient-poset. Similarly, the term ``semilattice''  will also be applied to structures $(S;\sqcup,\leq)$ where $\leq$ is a preorder on $S$ such the quotient-structure under the induced equivalence relation $\equiv$ is a ``real'' semilattice with the partial order induced by $\leq$ (thus, we avoid  precise but more complicated terms like ``pre-semilattice''). We call preorders (or pre-semilattices) $P,Q$ {\em equivalent} (in symbols, $P\simeq Q$) if their quotient-posets (resp., quotient-semilattices) are isomorphic. For subsets $A,B\subseteq S$ of a preorder $(S;\leq)$ we write $A\equiv B$ if every element of $A$ is equivalent to some element of $B$ and vice versa.
This terminology is especially convenient in the situation (which is typical below) when one operation $\sqcup$ on $S$ induces the pre-semilattice structure for several preorders on $S$.

We associate with any poset $Q$ the preorder $(Q^*;\leq^*)$ where $Q^*$ is the set of non-empty finite subsets of $Q$, and $S\leq^*R$ iff $\forall s\in S\exists r\in R(s\leq r)$. Let $Q^\sqcup$ be the quotient-poset of $(Q^*;\leq^*)$ and $\sqcup$ be the operation  of supremum in $S$ induced by the operation of  union in $Q^*$. Then $Q^\sqcup$ is a d-semilattice the join-irreducible elements of which coincide with the elements induced  by the singleton sets in $Q^*$ (the new smallest element $\bot$ corresponds to the empty subset of $Q$); thus, $(I(Q^\sqcup);\leq^*)\simeq Q$. Any element of $Q^\sqcup_\bot$ has a canonical decomposition. If $Q$ is well founded then so is also $Q^\sqcup$. The construction $Q\mapsto Q^\sqcup$ is a functor from the category of preorders to the category of semilattices (see \cite{s18} for additional details). We will use the following easy fact.

\begin{Proposition}\label{slattice}
Let $f:Q\to I(S)$ be a monotone function from a poset $Q$ to  the set of join-irreducible elements of a semilattice $S$. Then there is a unique semilattice homomorphism $f^\sqcup:Q^\sqcup\to S$ extending $f$. If $f$ is an embedding and $S$ is distributive then $f^\sqcup$ is an embedding.
\end{Proposition}

\subsection{Labeled trees and forests}\label{trees}

In this section we recall a notation system for levels of the FH of $k$-partitions introduced in \cite{s12}, with some additions and  new facts.

Let $\omega^*$ be the set of finite strings of natural numbers including the empty string $\varepsilon$.  By a {\em tree} we mean a nonempty initial segment of $(\omega^*;\sqsubseteq)$ where $\sqsubseteq$ is the prefix relation. A tree $T$ is {\em normal} if $\tau(i+1)\in T$ implies $\tau i\in T$. The {\em rank} of a finite tree is the length of a longest chain in the tree.

It is often useful to consider ``abstract'' trees (a special kind of posets) along with the ``concrete'' trees defined above. Obviously, any finite abstract tree is isomorphic to a normal tree above. In some definitions below we implicitly use this obvious relation between ``abstract'' and ``concrete'' trees. This slightly abuse notation but simplifies intuitive understanding of the definitions.

For any finite  tree $T$ and any $\tau\in T$, define the tree $T(\tau)=\{\sigma\in \omega^*\mid\tau\cdot\sigma\in T\}$ where $\cdot$ is the concatenation of strings (often omitted). Then any non-singleton finite tree $T$  is determined by  the singleton tree $\{\varepsilon\}$ and the finite  trees $T(i)$, $i\in T$ of lesser tree ranks than $T$, so $T=\{\varepsilon\}\cup\bigcup_{i\in T}T(i)$. We will use this representation  in the proofs by induction on ranks. By a {\em forest} we mean an initial segment of $(\omega^*\setminus\{\varepsilon\};\sqsubseteq)$. Note that there is a unique empty forest, and for any forest $F$ there is a unique tree $T$ with $F=T\setminus\{\varepsilon\}$. Any non-empty forest $F$ is isomorphic to the  disjoint union of trees $F(i)$, $i\in\omega\cap F$.

Next we recall notation related to iterated labeled forests from \cite{s12}. Let $(Q;\leq)$ be a preorder; abusing notation we often denote it just by $Q$. A {\em $Q$-tree}  is a pair $(T,t)$ consisting of a finite normal tree $T\subseteq\omega^*$,  and a labeling $t:T\to Q$. Let  ${\mathcal T}(Q)$ denote the set of all finite $Q$-trees.  The {\em $h$-preorder} $\leq_h$ on  ${\mathcal T}(Q)$ is defined as follows: $(T,t)\leq_h (V,v)$, if there is a monotone function $f:(T;\sqsubseteq)\to(V;\sqsubseteq)$ satisfying $\forall
\tau\in T(t(\tau))\leq v(f(\tau)))$. For any $q\in Q$, let $s(q)=(\{\varepsilon\},q)$ be the singleton tree labeled by $q$. Then $q\leq_Qr$ iff $s(q)\leq_hs(r)$. Note that $(T,t)\leq_hs(r)$ iff $t(\tau)\leq_Qr$ for all $\tau\in T$, $s(q)\leq_h(V,v)$ iff $q\leq_Qv(\sigma)$ for some $\sigma\in V$, and, if both $T,V$ are non-singleton then $T\leq_hV$ iff ($t(\varepsilon)\leq_Qv(\varepsilon)$ and $T(i)\leq_hV$ for all $i\in\omega\cap T$) or ($t(\varepsilon)\not\leq_Qv(\varepsilon)$ and $T\leq_hV(j)$ for some $j\in\omega\cap T$). This characterises the relation $\leq_h$ by induction on the tree rank. The construction $Q\mapsto{\mathcal T}(Q)$ may be considered as an endofunctor on the category of preorders (see \cite{s18} for details); then $s$ becomes a natural transformation from the identity endofunctor to ${\mathcal T}$. 

By Subsection \ref{notation}, $s$ induces the semilattice homomorphism $s^\sqcup$ from $Q^\sqcup$ into $\mathcal{T}(Q)^\sqcup$ (note that the latter semilattice is naturally isomorphic to $\mathcal{F}(Q)$); we often simplify  $s^\sqcup$ to $s$. The elements of $\mathcal{F}(Q)$ are naturally identified with the finite $Q$-labeled forests, and $I(\mathcal{F}(Q);\sqcup,\leq_h)\equiv_h\mathcal{T}(Q)$. Along with the operation $\sqcup$, there is a natural binary operation $\cdot$ on $\mathcal{F}(Q)$ defined as follows: $F\cdot G$ is the labeled forest obtained from $F$ by putting a copy of $G$ above every leaf of $F$. Clearly, $F\sqcup G\leq_hF\cdot G$ and $(F\cdot G)\cdot H=F\cdot( G\cdot H)$. For any $(F,t)\in\mathcal{T}(Q)^\sqcup$, let $r(F)=\bigsqcup\{t(\tau)\mid\tau\in F\}$. Then $r:\mathcal{F}(Q)\to Q^\sqcup$ is a semilattice homomorphism  such that $r(s(x))\equiv x$ for every $x\in Q$. The relation $\leq$ on $Q^\sqcup$ induces the relation $\leq'$ on $\mathcal{F}(Q)$ as follows: $F\leq'G$ iff $r(F)\leq r(G)$. Clearly, $r$ induces an equivalence of semilattices $(\mathcal{F}(Q);\sqcup,\leq')$ and $Q^\sqcup$, $I(\mathcal{F}(Q);\sqcup,\leq')\equiv's(Q)$, and the relations $\leq,\leq'$ coincide on $s(Q)$. 

By a {\em minimal $Q$-forest} we mean a  $Q$-forest not
$h$-equivalent  to a $Q$-forest of lesser cardinality. The next
characterization of the minimal $Q$-forests was obtained in \cite{s12}, Proposition 8.3.

\begin{Proposition}\label{minq}
\begin{enumerate}\itemsep-1mm
 \item Any singleton $Q$-forest is minimal.
 \item  A non-singleton $Q$-tree $(T,t)$ is minimal iff the $Q$-forest
$T\setminus\{\varepsilon\}$ is minimal, $\forall i\in\omega\cap T(t(i)\not\leq t(\varepsilon))$, and if $\omega\cap T=\{0\}$ then $t(0),t(\varepsilon)$ are incomparable.
 \item  A $Q$-forest having at least two $Q$-trees is minimal iff all its $Q$-trees are minimal and
pairwise incomparable under $\leq_h$.
 \end{enumerate}
  \end{Proposition}

Minimal forests are often useful to simplify proofs. E.g., if $F$ is minimal and consists of trees $F_0,\ldots,F_n$ then $F\equiv_hF_0\sqcup\cdots\sqcup F_n$ is automatically a canonical decomposition. It is easy to see that two minimal $h$-equivalent $Q$-forests are isomorphic, and that Proposition \ref{minq} induces a unary function $min$ on $\mathcal{F}(Q)$ such that $min(F)$ is a minimal $Q$-forest equivalent to $F$ (computing of such $min(F)$ may be called minimization of $F$). 

The function $min$ may be ``computed'' by induction on $|F|$ as follows: if $|F|=1$, set $min(F)=F$; if $F$ consists of trees $T_0,\ldots,T_n$, $n>0$, then delete all trees which are $h$-below some other tree, and minimize (by induction) all the remaining trees; if $F$ is a tree and $|F|>1$, then compute $F_1=t(\varepsilon)\cdot min(F\setminus\{\varepsilon\})$, delete from the resulting tree $F_1$ all nodes $\tau\neq\varepsilon$ with $\forall\sigma(\varepsilon\sqsubset\sigma\sqsubseteq\tau\to t(\sigma)\leq t(\varepsilon))$ to obtain the tree $F_2$, compute $F_3=t(\varepsilon)\cdot min(F_2\setminus\{\varepsilon\})$, and delete the root in $F_3$ whenever $\omega\cap F_3=\{0\}$ and $t(\varepsilon)\leq t(0)$.

Note that the  labels of $min(F)$ are among the labels of $F$, $min(F)$ is a tree whenever $F$ is a tree, and $F$ is join-irreducible iff $min(F)$ is a tree. If the structure $(Q;\leq)$ is computably presentable then the previous paragraph becomes a real algorithm for computing $min(F)$ from a given $F$.  

Define the sequence $\{\mathcal{T}_m(Q)\}_{m<\omega}$ of preorders by induction on $m$ as follows: $\mathcal{T}_0(Q)=Q$  and $\mathcal{T}_{m+1}(Q)=\mathcal{T}(\mathcal{T}_m(Q))$; the preorder on $\mathcal{T}_m(Q)$ is denoted by $\leq_h$  (for $m=0$ we identify $\leq_h$ with $\leq_Q$). For any $m<\omega$, we have the preorder embedding $s:\mathcal{T}_m(Q)\to\mathcal{T}_{m+1}(Q)$ which is sometimes convenient to denote $s_m$. Setting $\mathcal{F}_m(Q)=\mathcal{T}_m(Q)^\sqcup$, we obtain  sequences $\{\mathcal{F}_m(Q)\}$ of semilattices, and $s:\mathcal{F}_m(Q)\to\mathcal{F}_{m+1}(Q)$, $r:\mathcal{F}_{m+1}(Q)\to\mathcal{F}_{m}(Q)$ of homomorphisms such that $r(s(x))\equiv_h x$ for every $x\in \mathcal{F}_{m}(Q)$. For all $n\leq m$, we define a preorder $\leq^n$ on $\mathcal{F}_{m}(Q)$ by induction: $\leq^0$ is $\leq_h$ and, for $n<m$,  $F\leq^{n+1}G$ iff $r(F)\leq^nr(G)$. For all $p\leq m$, we define  a binary operation $\cdot^p$ on $\mathcal{F}_{m+1}(Q)$ by induction: $F\cdot^0G=F\cdot G$ and, for $p<m$, $F\cdot^{p+1}G=s(r(F)\cdot^pr(G))$. 

We collect some properties of the introduced objects in the following proposition where   $I(\mathcal{F}_{m}(Q);\sqcup,\leq^n)$ is abbreviated to  $I_n(\mathcal{F}_{m}(Q))$ and $s^n$ is the $n$-th iteration of $s$, i.e. $s^0$ is the identity and $s^{n+1}=s\circ s^n$. 

\begin{Proposition}\label{iterh}
\begin{enumerate}\itemsep-1mm
 \item For all $m$ and $n\leq m$ we have: $(\mathcal{F}_{m}(Q);\sqcup,\leq^n)$ is a d-semilattice,   $I_n(\mathcal{F}_{m}(Q))\equiv^0s^n(\mathcal{T}_{m-n}(Q))$, and the relations $\leq^0,\ldots,\leq^n$ coincide on $s^n(\mathcal{T}_{m-n}(Q))$.
 \item For all $p\leq m$ and $F,G,H\in\mathcal{F}_{m+1}(Q)$ we have: $r(F\cdot^0G)\equiv^0F\sqcup G$, $r(F\cdot^{p+1}G)\equiv^0F\cdot^p G$ for $p<m$, and $(F\cdot^p G)\cdot^p H\equiv^0F\cdot^p( G\cdot^p H)$.
 \end{enumerate}
\end{Proposition}

{\em Proof.} All  assertions follow by induction from the corresponding definitions and facts above, so we consider only the associativity of  $\cdot^p$, as an example. For $p=0$ the assertion is clear, so we assume it for $\cdot^p$, $p<m$, and check it for $\cdot^{p+1}$:
 \begin{eqnarray*}
  (F\cdot^{p+1} G)\cdot^{p+1} H=s(r(F)\cdot^{p}r(G))\cdot^{p+1} H=
 s(rs(r(F)\cdot^{p} r(G))\cdot^{p} r(H))\equiv^0\\  s((r(F)\cdot^{p} r(G))\cdot^{p} r(H))\equiv^0  s(r(F)\cdot^{p}( r(G)\cdot^{p} r(H)))\equiv^0\\ s(r(F)\cdot^{p}rs( r(G)\cdot^{p} r(H)))\equiv^0 s(r(F)\cdot^{p}r( G\cdot^{p+1} H))=F\cdot^{p+1} (G\cdot^{p+1} H).
 \end{eqnarray*}
 \qed


The sets $\mathcal{T}_m(Q)$, $m\geq0$, are pairwise disjoint but, identifying  $q\in Q$ with the corresponding singleton tree $s(q)$, we may think that $\mathcal{T}_0(Q)\subseteq\mathcal{T}_1(Q)$ and, moreover, $\mathcal{T}_0(Q)\subseteq\mathcal{T}_1(Q)\subseteq\cdots\subseteq\mathcal{T}_\omega(Q)$ and $\mathcal{T}_\omega(Q)$ has an induced preorder  also denoted by $\leq_h$. In a similar way we may think that $\mathcal{F}_\omega(Q)$ is a smallest structure containing any $\mathcal{F}_m(Q)$, $m<\omega$, as a substructure. This informal construction from some of my previous papers can be made precise (on the cost of additional technical details) by using the  known fact that the category of preorderes has arbitrary colimits. Since we need the precise construction in Section \ref{app}, we briefly recall some details of this construction.

Namely, we consider $\mathcal{T}_\omega(Q)$ as the colimit of the sequence of morphisms $g_m:\mathcal{T}_m(Q)\to\mathcal{T}_{m+1}(Q)$ defined by induction as follows: $g_0=s$, and $g_{m+1}(T,t)=(T,g_m\circ t)$. For any $T\in\mathcal{T}_\omega(Q)$, let $g(T)=g_m(T)$ where $m$ is the unique integer with $T\in\mathcal{T}_m(Q)$. Then $g$ is a function on $\mathcal{T}_\omega(Q)=\bigcup_{m}\mathcal{T}_m(Q)$ (where the summands are pairwise disjoint). Let $\leq$ be the smallest preorder on $\mathcal{T}_\omega(Q)$ such that $g(T)\equiv T$ for every $T\in\mathcal{T}_\omega(Q)$, and $T\leq V$ for all $T,V\in\mathcal{T}_m(Q)$ with $T\leq_hV$. Then the inclusions $\mathcal{T}_0(Q)\subseteq\mathcal{T}_1(Q)\subseteq\cdots\subseteq\mathcal{T}_\omega(Q)$ hold up to  $\equiv$ (intuitively, $\leq$ coincides with $\leq_h$ above). In checking this one has to use functions $g_{m,n}$, $m\leq n<\omega$, defined by $g_{m,n}=g_{n-1}\circ\cdots\circ g_{n}$ (for $m=n$, $g_{m,m}$ is the identity function on $\mathcal{T}_m(Q)$); thus,  $g_{m,n}$ is an embedding of $(\mathcal{T}_m(Q);\leq_h)$ into $(\mathcal{T}_n(Q);\leq_h)$.

A similar construction in the category of (pre-)semilattices yields the structure $\mathcal{F}_\omega(Q)$ as the colimit    of semilattice embeddings $g_{m,n}^\sqcup:(\mathcal{F}_m(Q);\sqcup,\leq_h)\to(\mathcal{F}_n(Q);\sqcup,\leq_h)$. The supremum operation $\sqcup$ on $\mathcal{F}_\omega(Q)$ is defined by $F\sqcup G=g_{m,p}(F)\sqcup_p g_{m,n}(G)$ where $p=max\{m,n\}$, $F\in\mathcal{F}_m(Q)$, and $G\in\mathcal{F}_n(Q)$. One easily checks that $F\leq F'$ and $G\leq G'$ imply $F\sqcup F'\leq G\sqcup G'$, hence  $(\mathcal{F}_\omega(Q);\sqcup,\leq)$ is a $d$-semilattice equivalent to $\mathcal{T}_\omega(Q)^\sqcup$ that contains $(\mathcal{F}_m(Q);\sqcup,\leq_h)$, $m<\omega$, as a substructure    and satisfies $I(\mathcal{F}_\omega(Q);\sqcup,\leq)\equiv\mathcal{T}_\omega(Q)$. 

The functions and relations defined before Proposition \ref{iterh} on  $\mathcal{T}_m(Q)$ and $\mathcal{F}_m(Q)$, $m<\omega$, are extended to  $\mathcal{T}_\omega(Q)$ and  $\mathcal{F}_\omega(Q)$ in the natural way. For any $T\in\mathcal{T}_\omega(Q)$, let $s(T)=s_m(T)$ where $m$ is the unique number with $T\in\mathcal{T}_m(Q)$. Then $s$ is a unary function on $\mathcal{T}_\omega(Q)$ that is monotone w.r.t. $\leq$ (because $g(s(T))\equiv s(g(T))$) and extends all $s_m$. It has the following properties: $s(q)\equiv q$ for $q\in Q$, $T<s(T)$ for every $T\in\mathcal{T}_\omega(Q)$ with $\forall q\in Q(T\not\equiv q)$, and $T\leq V$ iff $s(T)\leq s(V)$ for $T,V\in\mathcal{T}_\omega(Q)$. The function $s$ is extended to the semilattice embedding $s^\sqcup:\mathcal{F}_\omega(Q)\to\mathcal{F}_\omega(Q)$ which for simplicity is also  denoted by $s$. The binary operation $\cdot$ on $\mathcal{F}_\omega(Q)$ is defined by $F\cdot G=g_{m,p}(F)\cdot_p g_{m,n}(G)$ where $p=max\{m,n\}$, $F\in\mathcal{F}_m(Q)$, and $G\in\mathcal{F}_n(Q)$. One easily checks that $F\leq F'$ and $G\leq G'$ imply $F\cdot F'\leq G\cdot G'$, $F\sqcup G\leq F\cdot G$, and $(F\cdot G)\cdot H\equiv F\cdot( G\cdot H)$.

For any $F\in\mathcal{F}_m(Q)$, let $r(F)=r_m(F)$ for $m>0$ and $r(F)=F$ otherwise. Then $r$ is a unary function on $\mathcal{F}_\omega(Q)$ that is monotone w.r.t. $\leq$ (because $g(r(T))\equiv r(g(T))$) and extends all $r_m$. Furthermore, it is a semilattice epimorphism from $(\mathcal{F}_\omega(Q);\sqcup,\leq)$ onto itself such that $r(s(F))\equiv F$ for every $F\in\mathcal{F}_\omega(Q)$. Using this function $r$, we can define, for any $n<\omega$,  relations $\leq^n$ and operations $\cdot^n$ on $\mathcal{F}_\omega(Q)$ by induction as follows: $\leq^0$ is $\leq_h$,  $F\leq^{n+1}G$ iff $r(F)\leq^nr(G)$, $F\cdot^0G=F\cdot G$, and $F\cdot^{n+1}G=s(r(F)\cdot^nr(G))$. 
The next proposition  follows from Proposition \ref{iterh}. 

\begin{Proposition}\label{iterh1}
\begin{enumerate}\itemsep-1mm
 \item For any $n<\omega$ we have: $(\mathcal{F}_\omega(Q);\sqcup,\leq^n)$ is a d-semilattice, and the relations $\leq^0,\ldots,\leq^n$ coincide on  $s^n(\mathcal{T}_\omega(Q))\equiv^0I(\mathcal{F}_\omega(Q);\sqcup,\leq^n)$.
 \item For any $p<\omega$ and $F,G,H\in\mathcal{F}_\omega(Q)$ we have: $r(F\cdot^0G)\equiv^0F\sqcup G$, $r(F\cdot^{p+1}G)\equiv^0F\cdot^p G$, and $(F\cdot^p G)\cdot^p H\equiv^0F\cdot^p( G\cdot^p H)$.
 \end{enumerate}
\end{Proposition}
 

The preorder $Q$ is a {\em well quasiorder} (WQO) if it has neither infinite descending chains nor infinite antichains. A famous Kruskal's theorem implies that if $Q$ is WQO then $({\mathcal T}_Q;\leq_h)$ is WQO. It follows that any preorder $\mathcal{T}_m(Q)$, $m<\omega$, is WQO. It might be shown that the colimit preorder $(\mathcal{T}_\omega(Q);\leq)$ (hence also $(\mathcal{T}_\omega(Q);\leq^n)$ for every $n<\omega$) is  WQO. The preorders  $(\mathcal{F}_\omega(Q);\leq^n)$,  $n<\omega$, are also  WQOs.

Below we will mainly deal with the particular case $Q=\bar{k}$ of antichain of size  $2\leq k\leq\omega$. In this case $\mathcal{T}_m(\bar{k})\sqsubseteq\mathcal{T}_{m+1}(\bar{k})$ (i.e. the quotient-poset of $\mathcal{T}_m(\bar{k})$ is an initial segment of the quotient-poset of $\mathcal{T}_{m+1}(\bar{k})$). Note that $\mathcal{T}_\omega(\bar{k})$ is WQO for $k<\omega$, $\mathcal{T}_\omega(\bar{\omega})$ is well founded but not WQO, and  $\mathcal{T}_\omega(\bar{2})\sqsubseteq\mathcal{T}_\omega(\bar{3})\sqsubseteq\cdots\sqsubseteq\mathcal{T}_\omega(\bar{\omega})=\bigcup\{\mathcal{T}_\omega(\bar{k})\mid 2\leq k\leq\omega\}$. For $k=2$,  the quotient-poset of $(\mathcal{T}_\omega(\bar{2});\leq_h)$ has order type $\bar{2}\cdot\varepsilon_0$ (see  Proposition 8.28 in \cite{s12}).

Initial segments of $(\mathcal {F}_1(\bar{k});\leq_h)$ for $k=2,3$ are depicted below.\footnote{I thank Anton Zhukov for the help with making the pictures.}

\begin{center}
\includegraphics[scale=0.6]{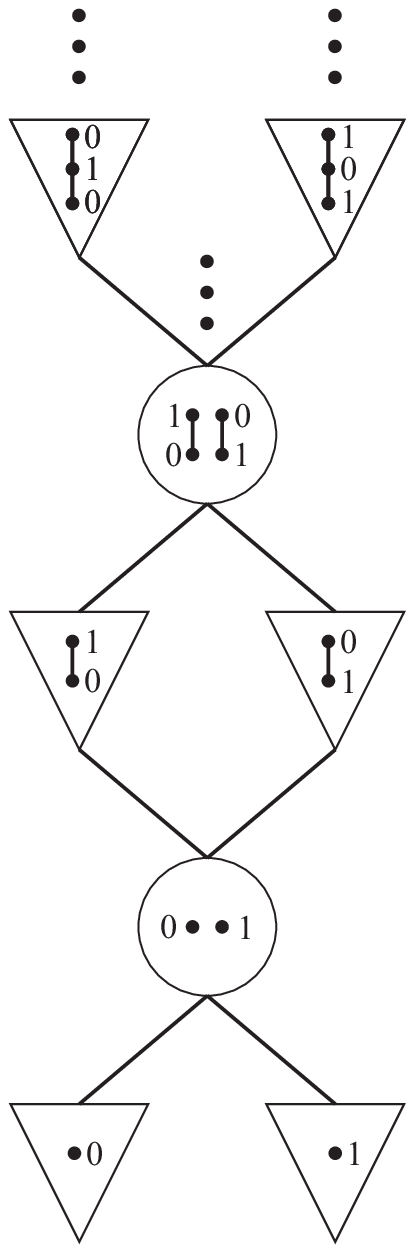}

Fig. 1. An initial segment of $(\mathcal {F}_1(\bar{2});\leq_h)$.
\end{center}

\begin{center}
\includegraphics[scale=0.5]{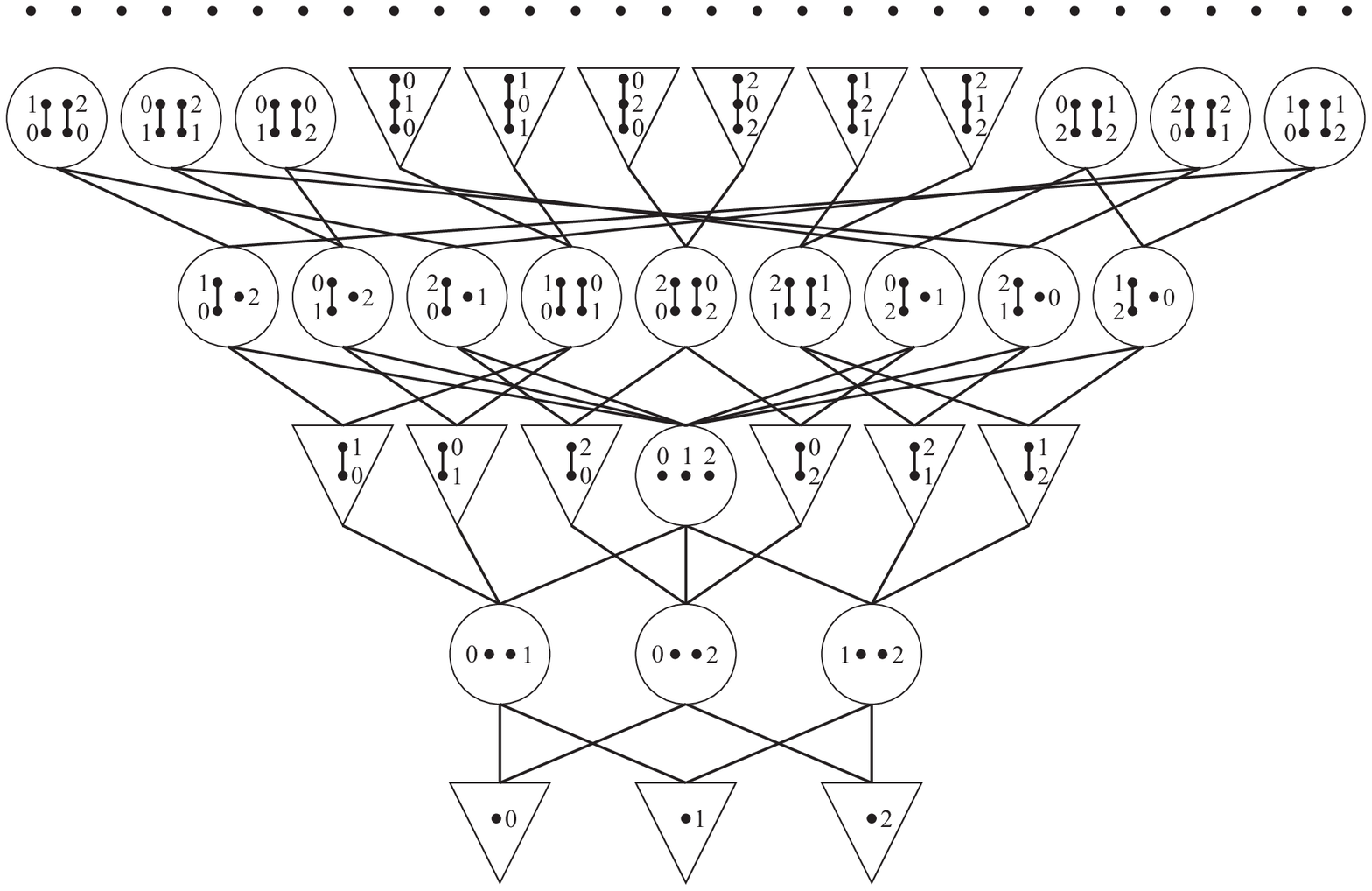}

Fig. 2. An initial segment of $(\mathcal{F}_1(\bar{3});\leq_h)$.
\end{center}

Note that while $\mathcal {F}_1(\bar{2})$ is semi-well-ordered with rank $\omega$, $\mathcal {F}_1(\bar{k})$ for $k>2$ is a wqo of rank $\omega$ having antichains of arbitrary finite size. The whole structure $\mathcal {F}_\omega(\bar{2})$  is also semi-well-ordered but with larger rank $\varepsilon_0$ (see  Proposition 8.28 in \cite{s12}). This structure is isomorphic to the FH of arithmetical sets in \cite{s83,s89} mentioned in the Introduction. The triangle levels (induced by trees) correspond to the ``non-self-dual'' $\Sigma$- and $\Pi$-levels of this hierarchy. More precisely, the $\Sigma$-levels (resp. $\Pi$-levels) correspond to (hereditary) 0-rooted (resp., 1-rooted) trees. According to Fig. 2,  the preorder $\mathcal{F}_1(\bar{k})$ for $k>2$ is much more complicated than for $k=2$. Nevertheless, the (generalised) non-self-dual levels of the corresponding FHs of $k$-partitions will again correspond to trees (depicted as triangles).

We conclude this subsection by defining minimal iterated $k$-forests $F\in\mathcal{F}_\omega(\bar{k})$. The definition is by induction on the unique $m$ with $F\in\mathcal{F}_m(\bar{k})$. For $m=0$, $F$ is minimal, if among its ``trees'' (which are of the form $i$, $i<k$) there are no repetitions. For $m>0$, $F$ is minimal, if it is a minimal $\mathcal{F}_{m-1}(\bar{k})$-forest in the sense of Proposition \ref{minq}, and all labels of $F$ are minimal iterated $k$-forests. As for the $Q$-forests, there is again a unary function $min$ on $\mathcal{F}_\omega(\bar{k})$ which minimizes iterated $k$-forests. It is computed in the obvious way: For $m=0$, $min(F)$ is obtained from $F$ by deleting repetitions of ``trees''; for $m>0$, $min(F)$ is obtained from $F$ by minimizing $F$ as a $\mathcal{F}_{m-1}(\bar{k})$-forest, and then by minimizing all labels of the obtained forest.

\section{Effective  Wadge hierarchy}\label{wkpart}

Since the EWH is a special case of the FH, we  first recall in this section some  information about the FH from \cite{s12,s19}, and then specialize it to obtain the EWH. 

We warn the reader that our definition of EWH uses set operations instead of the Wadge reducibility the reader could expect to see. The  Wadge reducibility (and especially the effective Wadge reducibility) leads to complex degree structures in non-zero-dimensional spaces which hide the hierarchy  (see \cite{s19a,s19} for more detailed discussions), so it cannot be used to define the (effective) WH in a broad enough context. The set-theoretic definition of the FH of sets was first introduced in \cite{s89} and then extensively studied in different contexts (see \cite{s08a} for a partial survey). The extension of the FH to $k$-partitions was first introduced in \cite{s12}.

By a {\em base in a set $X$} we mean a sequence $\mathcal{L}(X)=\{\mathcal{L}_n\}_{n<\omega}$ of subsets of $P(X)$ such that any $\mathcal{L}_n$ is closed under union and intersection, contains $\emptyset,X$, and  $A\in\mathcal{L}_n$ implies that $A,X\setminus A\in\mathcal{L}_{n+1}$. The base $\mathcal{L}(X)$ is {\em reducible} if  every its level $\mathcal{L}_n$ has the reduction property.

With any base $\mathcal{L}(X)$  we  associate some other bases as follows. For any $m<\omega$, let $\mathcal{L}^m(X)=\{\mathcal{L}_{m+n}(X)\}_n$; we call this base   {\em $m$-shift of $\mathcal{L}(X)$}.  For any $U\in\mathcal{L}_0$, let $\mathcal{L}(U)=\{\mathcal{L}_n(U)\}_{n<\omega}$ where $\mathcal{L}_n(U)=\{U\cap S\mid S\in\mathcal{L}_{n}(X)\}$; we call this base in $U$ the {\em $U$-restriction of $\mathcal{L}(X)$}.  

We define the FH not only of subsets of $X$ but also  of $k$-partitions  $A:X\to \bar{k}$, $1<k<\omega$. Note that $2$-partitions of $X$ are essentially subsets of $X$. 
For any finite tree $T\subseteq\str$ and any $T$-family $\{U_\tau\}$ of subsets of $X$, we define the $T$-family $\{\tilde{U}_\tau\}$ of subsets of $X$ by $\tilde{U}_\tau=U_\tau\setminus\bigcup\{U_{\tau'}\mid\tau\sqsubset\tau'\in T\}$. The $T$-family $\{U_\tau\}$ is {\em monotone} if $U_\tau\supseteq U_{\tau'}$ for all $\tau\sqsubseteq\tau'\in T$. We associate with any $T$-family $\{U_\tau\}$ the monotone $T$-family $\{U'_\tau\}$ by $U'_\tau=\bigcup_{\tau'\sqsupseteq\tau}U_{\tau'}$. A $T$-family $\{V_\tau\}$  is {\em reduced} if it is monotone and satisfies $V_{\tau i}\cap V_{\tau j}=\emptyset$ for all $\tau i,\tau j\in T$. Obviously, for any reduced $T$-family $\{V_\tau\}$  the components $\tilde{V}_\tau$ are pairwise disjoint. 

We will use the following technical notions. The first one is the notion ``$F$ is a  $T$-family in $\mathcal{L}(X)$''  defined by induction as follows: if $T\in\mathcal{T}_0(\bar{k})$ then $F=\{X\}$; if $(T,t)\in\mathcal{T}_{m+1}(\bar{k})$ then $F=(\{U_\tau\},\{F_\tau\})$ where $\{U_\tau\}$ is a monotone $T$-family of  ${\mathcal L}_0$-sets with $T_\varepsilon=X$ and, for each $\tau\in T$, $F_\tau$ is  a  $t(\tau)$-family in $\mathcal{L}^1(\tilde{U}_\tau)$. The version of this notion ``$F$ is a reduced $T$-family in $\mathcal{L}(X)$'' is obtained by taking the reducible $T$-families in place of the monotone ones.
The second is the notion  ``a  $T$-family $F$ in $\mathcal{L}(X)$ determines $A:X\to\bar{k}$'' defined by induction as follows:
if $T\in\mathcal{T}_0(\bar{k})$, $T=i<k$ (so $F=\{X\}$), then $T$ determines the constant partition $A=\lambda x.i$;
if $(T,t)\in\mathcal{T}_{m+1}(\bar{k})$ (so $F$ is of the form $(\{U_\tau\},\{F_\tau\})$) then $T$ determines the $k$-partition $A$ such that $A|_{\tilde{U}_\tau}=B_\tau$  for every $\tau\in T$, where $B_\tau:\tilde{U}_\tau\to\bar{k}$ is the $k$-partition of $\tilde{U}_\tau$ determined by $F_\tau$. 

As explained in \cite{s19a}, the $T$-family $F$  that determines $A$ provides a mind-change algorithm for computing $A(x)$ (see Section 3 of \cite{s19} for additional details).
We are ready to give a precise definition of  the FH of $k$-partitions over ${\mathcal L}(X)$.

\begin{Definition}\label{fh}
The FH of $k$-partitions over ${\mathcal L}(X)$ is the family $\{\mathcal{L}(X,T)\}_{T\in\mathcal{T}_\omega(\bar{k})}$ of subsets of $k^X$ where ${\mathcal L}(X,T)$ is the set of $A:X\to\bar{k}$ determined by some $T$-family in $\mathcal{L}(X)$. Let red-${\mathcal L}(X,T)$ denote the set of $A:X\to\bar{k}$ determined by some reduced $T$-family in $\mathcal{L}(X)$. 
\end{Definition}

As shown in \cite{s12}, $T\leq_hS$ implies $\mathcal{L}(X,T)\subseteq\mathcal{L}(X,S)$, hence $(\{\mathcal{L}(X,T)\mid T\in\mathcal{T}_\omega(\bar{k})\};\subseteq)$ is WQO. 
The FH of sets  obtained from this construction for $k=2$ is even semi-well-ordered since the quotient-poset of $(\mathcal{T}_2(\omega);\leq_h)$ has order type $\bar{2}\cdot\varepsilon_0$ (see Definition 8.27 and Proposition 8.28 in \cite{s12}). 

The FH of $k$-partitions over the effective Borel base $\mathcal{L}(X)=\{\Sigma^0_{1+n}(X)\}$ in an effective cb$_0$-space $X$ is written as  $\{\Sigma(X,T)\}_{T\in\mathcal{T}_\omega(\bar{k})}$ and called the {\em effective Wadge hierarchy in $X$}. For $k=2$ the structure of levels degenerate to the semi-well-ordered structure  which enables the $\Sigma, \Pi$-notation for them. The EWH of sets subsumes many hierarchies including those mentioned in Section \ref{prel}.

The corresponding boldface FH $\{\mathbf{\Sigma}(X,T)\}_{T\in\mathcal{T}_\omega(\bar{k})}$ over the finite Borel base $\mathcal{L}(X)=\{\bfSig^0_{1+n}(X)\}$ is written as $\{\bfSig(X,T)\}_{T\in\mathcal{T}_\omega(\bar{k})}$ and is called  {\em finitary Wadge hierarchy in $X$}. It forms a small but important fragment of the whole (infinitary) Wadge hierarchy of $k$-partitions in $X$ (which is constructed from the whole Borel hierarchy by taking countable well-founded trees $T$ in place of the finite trees $T$). The latter hierarchy, which may be defined  in arbitrary space, was introduced and studied in \cite{s19a}. For the effective and boldface versions we have the obvious inclusions $\Sigma(X,T)\subseteq\mathbf{\Sigma}(X,T)$. In this paper we stick to levels of EWH corresponding to finite trees; the levels corresponding to computable well-founded trees (briefly discussed in \cite{s19}) are important on their on and we plan to investigate them in a separate publication.

We collect some simple properties of the EWH which are either contained in \cite{s12,s19} or easily follow from the definitions. As above, the effective Borel base $\{\Sigma^0_{1+n}(X)\}$ in $X$ is sometimes abbreviated to ${\mathcal L}$.

\begin{Proposition}\label{fineprop}
\begin{enumerate}\itemsep-1mm
\item $T\leq_hV$ implies $\Sigma(X,T)\subseteq\Sigma(X,V)$, hence $(\{\Sigma(X,T)\mid T\in\mathcal{T}_\omega(\bar{k})\};\subseteq)$ is WQO for $k<\omega$ and is well founded for $k=\omega$.
\item If  ${\mathcal L}$ is reducible then $\mathcal{L}(X,T)=$red-$\mathcal{L}(X,T)$ for every $T\in\mathcal{T}_\omega(\bar{k})$.
\item If $(T,t)\in\mathcal{T}_{m+1}(\bar{k})$ is a non-singleton normal tree, $\omega\cap T=\{i\mid i<p\}$, and $A\in\mathcal{L}^n(X,T)$ is determined by  an $\mathcal{L}^n$-family $(\{U_{\tau_0}\},\{U_{\tau_0\tau_1}\},\ldots)$ then $A|_{\overline{U}_\varepsilon}\in\mathcal{L}^{n+1}(\overline{U}_\varepsilon,t(\varepsilon))$ and $A|_{U_i}\in\mathcal{L}^n(U_i,T(i))$ for each $i<p$.
\item Let $A\in k^X$, $T$ be as in the previous item and $V_0,\ldots,V_{p-1}$ be pairwise disjoint $\mathcal{L}_n$-sets such that $A|_{V_i}\in\mathcal{L}^n(V_i,T(i))$ for $i<p$, and $A|_{\overline{V}}\in\mathcal{L}^{n+1}(\overline{V},t(\varepsilon))$ where $V=V_0\cup\cdots\cup V_{p-1}$. Then $A\in\mathcal{L}^n(X,T)$. 
\item Any level $\Sigma(X,T)$  is closed downwards under the effective Wadge reducibility  on $k^X$.
\end{enumerate}
\end{Proposition}

In \cite{s19a,s19} the following preservation property for levels of the (effective) WH was established.

\begin{Proposition}\label{pres}
Let $f:Y\to X$ be a computable effectively open surjection between  effective cb$_0$-spaces. Then, for all $T\in\mathcal{T}_\omega(\bar{k})$ and $A\in k^X$, we have: $A\in\Sigma(X,T)$ iff $A\circ f\in\Sigma(Y,T)$. Similarly for the boldface versions and continuous open surjections between cb$_0$-spaces. 
\end{Proposition}

\section{Non-collapse property}\label{ncol}

Here we establish some general facts about the non-collapse property. First we carefully define  natural versions of this property.

 We say that  EWH $\{\Sigma(X,T)\}_{T\in\mathcal{T}_\omega(\bar{k})}$ {\em does not collapse at level $T$} if $\Sigma(X,T)\not\subseteq\Sigma(X,V)$ for each $V\in\mathcal{T}_\omega(\bar{k})$ with $T\not\leq_hV$; the hierarchy  {\em strongly does not collapse at level $T$} if $\Sigma(X,T)\not\subseteq\bigcup\{\Sigma(X,V) \mid V\in\mathcal{T}_\omega(\bar{k}),\;T\not\leq_hV\}$. We say that $\{\Sigma(X,T)\}_{T\in\mathcal{T}_\omega(\bar{k})}$ {\em (strongly) does not collapse} if it (strongly) does not collapse at any level $T\in\mathcal{T}_\omega(\bar{k})$.   
 
 Note that for the case of sets $k=2$ these definitions are equivalent to the standard definition of non-collapse in DST ($\Sigma$-levels are distinct from the corresponding $\Pi$-levels), and the strong version is equivalent to the non-strong one. Thus, the strong versions are not visible for $k=2$.
 
\begin{Proposition}\label{ncolprop}
\begin{enumerate}
 \item The EWH $\{\Sigma(X,T)\}$  does not collapse iff the quotient-poset of $(\mathcal{T}_\omega(\bar{k});\leq_h)$ is isomorphic to  $(\{\Sigma(X,T)\mid T\in\mathcal{T}_\omega(\bar{k})\};\subseteq)$.
 \item If the EWH in $X$  does not collapse at level $T$ and $\Sigma(X,T)$ has a complete element w.r.t. the effective Wadge reducibility then the hierarchy  strongly does not collapse at level $T$.
 \item If the EWH in $X$  does not collapse and every its level has a complete element w.r.t. the effective Wadge reducibility then the hierarchy  strongly does not collapse.
 \end{enumerate}
 \end{Proposition}
 
 {\em Proof.} (1) Obvious from the definition.
 
(2) Let $A$ be effective Wadge complete in $\Sigma(X,T)$ and let $V\in\mathcal{T}_\omega(\bar{k})$ satisfy $T\not\leq_hV$; it suffices to show that $A\not\in\Sigma(X,V)$. Since the hierarchy does not collapse at level $T$, there is $B\in\Sigma(X,T)\setminus\Sigma(X,V)$. Since $A$ is complete in $\Sigma(X,T)$, we have $B\leq_{eW}^XA$. Since $\Sigma(X,V)$ is closed downwards under $\leq_{eW}^X$ by Proposition \ref{fineprop}(5), we get $A\not\in\Sigma(X,V)$.

(3) Follows from (2).
 \qed

The non-collapse for the boldface versions are defined in the same way, and the analogue of the previous proposition holds for them with essentially the same proof, including the version for infinitary FH. 

In the effective case, there are also the following uniform versions of  non-collapse property which relate EWH to the corresponding WH. The EWH $\{\Sigma(X,T)\}$ {\em  uniformly does not collapse at level $T$} if $\Sigma(X,T)\not\subseteq\mathbf{\Sigma}(X,V)$ for each $V\in\mathcal{T}_\omega(\bar{k})$ with $T\not\leq_hV$. The EWH $\{\Sigma(X,T)\}$ {\em  uniformly does not collapse}  if it uniformly does not collapse at every level. The hierarchy   {\em strongly uniformly does not collapse} at level $T$ if $\Sigma(X,T)\not\subseteq\bigcup\{\mathbf{\Sigma}(X,V) \mid V\in\mathcal{T}_\omega(\bar{k}),\;T\not\leq_hV\}$. The hierarchy   {\em strongly uniformly does not collapse} if $\Sigma(X,T)\not\subseteq\bigcup\{\mathbf{\Sigma}(X,V) \mid V\in\mathcal{T}_\omega(\bar{k}),\;T\not\leq_hV\}$ for all $T\in\mathcal{T}_\omega(\bar{k})$.

The following analogue of Proposition \ref{ncolprop} holds for the uniform version, essentially with the same proof.

\begin{Proposition}\label{ncolprop1}
\begin{enumerate}
 \item If the EWH in $X$ uniformly does not collapse at level $T$ and  $\Sigma(X,T)$ has a complete element w.r.t. $\leq_{eW}^X$ which is also Wadge complete in $\mathbf{\Sigma}(X,T)$ then the hierarchy strongly uniformly does not collapse at level $T$.
 \item If the EWH in $X$  uniformly does not collapse and every its level has an element  with the properties described in (1) then the hierarchy strongly uniformly does not collapse.
 \end{enumerate}
 \end{Proposition}
 
The next assertion follows from the inclusions between levels.

\begin{Proposition}\label{ncolinc}
For any effective cb$_0$-space $X$ we have: if $\{\Sigma(X,T)\}$ (strongly) uniformly does not collapse (at level $T$) then both $\{\Sigma(X,T)\}$ and $\{\mathbf{\Sigma}(X,T)\}$ (strongly) do not collapse (at level $T$).
\end{Proposition}

For cb$_0$-spaces $X$ and $Y$, let $X\leq_{co}Y$ mean that there is a continuous open surjection $f$ from $Y$ onto $X$. For effective cb$_0$-spaces $X$ and $Y$, let $X\leq_{eco}Y$ mean that there is a computable effectively open surjection $f$ from $Y$ onto $X$. Clearly, both $\leq_{eco}$ and $\leq_{co}$ are preorders, and the first preorder is contained in the second. The non-collapse property is inherited w.r.t. these preorders:

\begin{Proposition}\label{ncolpres}
\begin{enumerate}
 \item If $X\leq_{co}Y$ and $\{\mathbf{\Sigma}(X,T)\}_{T\in\mathcal{T}_\omega(\bar{k})}$ (strongly) does not collapse (at level $T$) then $\{\mathbf{\Sigma}(Y,T)\}$ (strongly) does not collapse (at level $T$). The same holds for the  infinitary version of  WH in $X$.
 \item If $X\leq_{eco}Y$ and $\{\Sigma(X,T)\}_{T\in\mathcal{T}_\omega(\bar{k})}$ (strongly)  does not collapse (at level $T$) then $\{\Sigma(Y,T)\}$ (strongly) does not collapse (at level $T$). The same holds for the uniform version of non-collapse property.
 \end{enumerate}
 \end{Proposition}
 
 {\em Proof.} All assertions  follow from the definitions and  the preservation property, so consider only  the finitary version in item (1). Let $X\leq_{co}Y$ via $f:Y\to X$, and $\{\mathbf{\Sigma}(X,T)\}$  does not collapse at level $T$. We have to show that $\mathbf{\Sigma}(Y,T)\not\subseteq\mathbf{\Sigma}(Y,V)$ for any fixed $V\in\mathcal{T}_\omega(\bar{k})$ with $T\not\leq_hV$. Choose $A\in\mathbf{\Sigma}(X,T)\setminus\mathbf{\Sigma}(X,V)$. By Proposition \ref{pres} we get $A\circ f\in\mathbf{\Sigma}(Y,T)\setminus\mathbf{\Sigma}(Y,V)$.
 \qed

\begin{Corollary}\label{ncolpres1}
\begin{enumerate}
 \item If $X$ is quasi-Polish and $\{\mathbf{\Sigma}(X,T)\}_{T\in\mathcal{T}_\omega(\bar{k})}$ (strongly) does not collapse (at level $T$) then $\{\mathbf{\Sigma}(\calN,T)\}$  (strongly)  does not collapse (at level $T$). The same holds for the infinitary version of the WH.
 \item If $X$ is computable quasi-Polish and $\{\Sigma(X,T)\}_{T\in\mathcal{T}_\omega(\bar{k})}$ (strongly)  does not collapse (at level $T$) then $\{\Sigma(\calN,T)\}$ (strongly) does not collapse (at level $T$). The same holds for the uniform version.
  \item  If $X$ is the product of a sequence $\{X_n\}$ of nonempty cb$_0$-spaces, and the finitary WH $\{\mathbf{\Sigma}(X_n,T)\}$  (strongly)  does not collapse (at level $T$) for some $n<\omega$, then $\{\mathbf{\Sigma}(X,T)\}$ (strongly)   does not collapse (at level $T$). The same holds for the infinitary version of the WH.
  \item If $X$ is the product of a uniform sequence $\{X_n\}$ of nonempty effective cb$_0$-spaces, and $\{\Sigma(X_n,T)\}$ (strongly)  does not collapse (at level $T$) for some $n<\omega$, then $\{\Sigma(X,T)\}$  (strongly) does not  collapse (at level $T$). The same holds for the uniform version. 
 \end{enumerate}
 \end{Corollary}
 
 {\em Proof.} (1) Follows from Proposition \ref{ncolpres}(1) since $X$ is quasi-Polish iff $X\leq_{co}\calN$.
 
 (2) Follows from Proposition \ref{ncolpres}(2) since $X$ is computable quasi-Polish iff $X\leq_{eco}\calN$.
 
 (3) Follows from Proposition \ref{ncolpres}(1) since  $X_n\leq_{co}X$.
 
 (4) Follows from Proposition \ref{ncolpres}(2) since  $X_n\leq_{eco}X$.
 \qed
 
Although the assertion  (1)  is void (because the infinitary WH in $\calN$  strongly  does not collapse \cite{wad84,km19}), it is of some methodological interest because it shows that proving the non-collapse of WH in any quasi-Polish space is at least as complicated as proving it in $\calN$, and the proof of the latter fact is highly non-trivial. The same applies to item (2) but this assertion is non-void because the non-collapse of EWH in $\calN$ was open until this paper, to my knowledge. 

In the remaining sections we give prominent examples of spaces with the non-collapse property. A good strategy to obtain broad classes of such spaces is to make them as low as possible w.r.t. $\leq_{co},\leq_{eco}$, and use the preservation property.

\section{Effective Wadge hierarchy in $\mathbb{N}$}\label{app}

In this section we discuss the EWH of $k$-partitions in $\mathbb{N}$, for any fixed number $k$ satisfying $2\leq k\leq\omega$. This hierarchy coincides with the FH of arithmetical $k$-partitions of $\omega$ first considered in \cite{s83}. In particular, we prove that the EWH of $k$-partitions in $\mathbb{N}$ strongly does not collapse.

\subsection{Complete numberings}\label{precomp}

Here we recall some information on precomplete and complete numberings from \cite{ma70,er77}, their strong relativizations introduced in \cite{s83}, and prove some additional facts. 
 
First we fix some notation and recall some notions of computability theory (we refer the reader to \cite{ro67,ma70,er77} for additional details). In particular, $\leq_T$ denotes the Turing reducibility on $P(\omega)$ (and also on $k^\omega$ for each $2\leq k<\omega$), $\leq_m$ denotes the many-one reducibility on $P(\omega)$, and $\equiv_T$, $\equiv_m$ denote the corresponding equivalence relations. As usual, by $\{\Sigma^0_{1+n}\}_{n<\omega}$ we denote the arithmetical hierarchy of subsets of $\omega$. As is well known, every level $\Sigma^0_{1+n}$ has the reduction property.

If the contrary is not specified explicitly, by a (partial) function we mean a (partial) function from $\omega$ to $\omega$, and by a set we mean a subset of $\omega$. For a partial function $\psi$, by
$\psi(x)\downarrow$ (resp. by $\psi(x)\uparrow$) we denote the fact that $x\in dom(\psi)$ (resp. $x\not\in dom(\psi)$). It is sometimes useful to identify a partial function $\psi$ with the total function $\psi:
\omega\rightarrow\omega\cup\{\perp\}$, where $\bot$ is a new element such that
$\psi(x)=\bot$ iff $\psi(x)\uparrow$.

By {\em a numbering} we mean any function $\nu$ with $dom(\nu)=\omega$, and by
{\em a numbering of a  set $S$} --- a numbering $\nu$ with $rng(\nu)=S$. A numbering $\mu$ is {\em reducible} to a numbering $\nu$ (in symbols, $\mu\leq\nu$) if there is a computable total function $f$ with $\mu=\nu\circ f$. Note that $\leq$ and $\leq_m$ coincide on $2^\omega$. For any $h\in\omega^\omega$,  $\leq^h$  denotes the $h$-relativization of the reducibility relation. The induced equivalence relations are denoted by $\equiv$ and $\equiv^h$, respectively.

Associate with any numberings $\mu,\nu$ and $\nu_k$ $(k<\omega)$ the numberings
$\mu\oplus\nu$,
$\mu\otimes\nu$ and $\bigoplus_k\nu_k$, called respectively  {\em the join of $\mu$
and $\nu$}, {\em the product of $\mu$ and $\nu$}, and {\em the infinite join of
$\nu_k$ $(k<\omega)$}, and defined as follows:
 $$(\mu\oplus\nu)2n=\mu n,\;(\mu\oplus\nu)(2n+1)=\nu n,\;
(\mu\otimes\nu) \langle
x,y\rangle=(\mu x,\nu y),\; (\bigoplus_k\nu_k)\langle x,y\rangle=\nu_x(y).$$

Here $\langle .,.\rangle$ is a standard computable bijection between
$\omega\times \omega$ and $\omega$; let $l,r$ be total computable functions such that $n=\langle 
 l(n),r(n)\rangle$ for every $n<\omega$. Similar notation is used for  coding  longer tuples.
The operations $\oplus,\otimes,\bigoplus$ are applicable to sets (by identifying sets with characteristic
functions considered as numberings). Iteration of $\oplus$ may be of course used to define join $\nu=\nu_0\oplus\cdots\oplus\nu_{p-1}$ of every finite sequence of numberings $\nu_0,\ldots,\nu_{p-1}=\bigoplus_{i<p}\nu_p$, $p>0$. We will use the following technically convenient syntactic version of this join:
$\nu(x)=\nu_i(a)$ where $a$ is the quotient and $i=rest(x,p)$ is the remainder of the division of $x$ by $p$, i.e. $x=p\cdot a+i$.

Following A.I. Maltsev \cite{ma70}, by $\varkappa$  we denote the Kleene numbering of the computable partial (c.p.) functions. More standard notation for $\varkappa$ is of course  $\varphi$ \cite{ro67}; we chose the less standard notation because it is used in the paper \cite{s83} which is often cited below. Note that $\tilde{\varkappa}=\bigoplus_k\varkappa_k$ is the standard universal c.p. function. By $\varkappa^h$ we denote the $h$-relativization of $\varkappa$ to any oracle $h\in\omega^\omega$, and by $h'=dom(\varkappa^h)$ --- the Turing jump of $h$ which may be considered as an operator on $k^\omega$ for every $2\leq k\leq\omega$. 

 
For any oracle $h\in\omega^\omega$, a numbering $\nu$   is called {\em $h$-precomplete}, if  for any $h$-computable partial ($h$-c.p.) function $\psi$   there is a total computable function  $t$
(called  a $\nu$-totalizer of   $\psi$) such that   $\nu t(x) =\nu\psi(x)$ for $x\in
dom(\psi)$. For $h=\emptyset$, $h$-precomplete numberings are called precomplete. 
Note  that  if  $\nu$   is  precomplete  and  $\{\psi_n\}$  is a uniform
sequence of c.p. functions then there is a uniform sequence $\{t_n\}$
of  $\nu$-totalizers for them. 
The precomplete numberings are just the numberings satisfying an effective version of the Kleene recursion theorem. 
A numbering  $\nu$  is {\em $h$-complete} (w.r.t.  $a\in rng(\nu)$), if for every $h$-c.p.  function $\psi$ there is  a computable function $t$ (called a $\nu$-totalizer of $\psi$
w.r.t. $a$) such that $\nu t(x) =   \nu\psi (x)$ for $\psi(x)\downarrow$
and   $\nu t(x) = a$ for $\psi(x)\uparrow$. It is known that if $\mu,\nu$ are $h$-precomplete (resp. $h$-complete) then so is also $\mu\otimes\nu$. Note that the constant numberings $\lambda n.s$, $s\in S$, are precisely the numberings in $S^\omega$ which are $h$-complete for every oracle $h$. Let $S^\omega_\bot$ be obtained from $S^\omega$ by adjoining a new element $\bot$ such that $\bot<^h\nu$ for all $h\in\omega^\omega$ and $\nu\in S^\omega$.

\begin{Proposition}\label{hprecomp}
\begin{enumerate}
\item The structure $(S_\bot^\omega;\oplus,\leq^h)$ is a distributive semilattice.

\item Any $h$-precomplete numbering $\nu:\omega\to S$ is join-irreducible in $(S_\bot^\omega;\oplus,\leq^h)$.

\item For any $h$-precomplete numberings $\mu,\nu$ we have: $\mu\leq^h\nu$ iff $\mu\leq\nu$.

\item For any finite families $\{\mu_i\}$ and $\{\nu_j\}$ of $h$-precomplete numberings we have: $\bigoplus_i\mu_i\leq^h\bigoplus_j\nu_j$ iff $\bigoplus_i\mu_i\leq\bigoplus_j\nu_j$.
 \end{enumerate}
\end{Proposition}

{\em Proof.} Items (1) and (2) are relativizations of known facts \cite{er77}, item (3) is clear from the definition. For item (4), using (2) and (3), we obtain: $\bigoplus_i\mu_i\leq^h\bigoplus_j\nu_j$ iff $\forall i(\mu_i\leq^h\bigoplus_j\nu_j)$ iff $\forall i\exists j(\mu_i\leq^h\nu_j)$ iff $\forall i\exists j(\mu_i\leq \nu_j)$ iff $\bigoplus_i\mu_i\leq\bigoplus_j\nu_j$.
 \qed

The Kleene numbering $\varkappa$ and its relativization $\varkappa^h$ are the most important for this section; $\varkappa^h$ is $h$-complete w.r.t. the empty function, and the universal $h$-c.p. function $\tilde{\varkappa}^h:\omega\to\omega_\bot$ is $h$-complete w.r.t. $\bot$. The Kleene recursion theorem for $\varkappa\otimes\varkappa$ is known as the double recursion theorem: for any uniformly computable sequences $\{g_x\}$, $\{h_x\}$ of partial functions there is $e$ such that $g_e=\varkappa_{l(e)}$ and $h_e=\varkappa_{r(e)}$.

For any $s\in S$ and $h\in\omega^\omega$, we defined in \cite{s83} the unary operation $p^h_s$ on $S^\omega$, which modifies ans relativises the completion operation from \cite{er77}, as follows: $[p^h_s(\nu)](x)=s$ if $\tilde{\varkappa}^h(x)\uparrow$ and $[p^h_s(\nu)](x)=\nu\tilde{\varkappa}^h(x)$ otherwise. Next we recall  an operation which, for $k=2$, is a natural combination of these operations and Turing jump.

\begin{Definition}\label{operG}
For \( \mu ,\nu \in S^\omega \) and \( h\in \omega^\omega  \),
let \( G(\mu ,\nu ,h)=\bigoplus _{n<\omega }p_{\nu (n)}^{h}(\mu ) \). In other words, we have: 
\( [G(\mu ,\nu ,h)]\langle n,x\rangle =
\nu (n) \) for \( \tilde{\varkappa} ^{h}(x)\uparrow  \) and
\( [G(\mu ,\nu ,h)]\langle n,x\rangle =
\mu \tilde{\varkappa} ^{h}(x) \) for \( \tilde{\varkappa} ^{h}(x)\downarrow  \). 
\end{Definition}

The function \( G:S^{\omega }\times S^{\omega }\times\omega^\omega
\to S^{\omega } \),  introduced in \cite{s83}, is crucial for this section. Its relations to the EWH in $\mathbb{N}$ is explained by the fact that, as also levels of this hierarchy, the operation is defined by using a relativized mind-change construct uniformly on oracles.

The next proposition collects some  facts from Propositions 1 --- 3 in \cite{s83}. Item (7), formulated here in a slightly modified form,  obviously remains true. 

\newcommand{\rrr}{G(\mu ,\nu ,h)}
\begin{Proposition}\label{p-r}
\begin{enumerate}
\item \( G(\mu ,\bigoplus _{n<\omega }\nu _{n},h)\equiv
\bigoplus _{n}G(\mu ,\nu _{n},h) \).

\item \( \mu \leq G(\mu ,\nu ,h) \), \( \nu \leq G(\mu ,\nu ,h) \),
\( G(\mu ,\nu ,h)\leq ^{h'}\mu \oplus \nu  \).


\item If \( \mu \leq ^{h}\mu _{1} \) then \( G(\mu ,\nu ,h)\leq G(\mu _{1},\nu ,h) \).

\item If \( \mu \leq ^{g}\mu _{1} \) and \( h^{\prime }\leq _{T}g \)
then \( G(\mu ,\nu ,h)\leq ^{g}G(\mu _{1},\nu ,h) \).

\item If \( \nu \leq ^{g}\nu _{1} \) then \( G(\mu ,\nu ,h)\leq^g G(\mu ,\nu _{1},h) \).


\item If \( h\leq _{T}h_{1} \) then \( G(\mu ,\nu ,h)\leq G(\mu ,\nu ,h_{1}) \).


\item For all \( f,g,h\in\omega^\omega \) we have: if $f=g=\lambda n.m$ for some $m<\omega$ then $G(f,g,h)=\lambda n.m$, otherwise
\( G(f,g,h)\equiv _{T}f\oplus g\oplus h^{\prime } \).

\item If \( \mu  \) is \( h^{\prime } \)-complete and \( \nu \leq ^{h^{\prime }}\mu  \)
then \( G(\mu ,\nu ,h)\equiv \mu  \).

\item If \( \nu  \) is \( h \)-complete w.r.t. \( s\in S \) then so is \( \rrr  \).

\item \( G(\mu ,G(\nu ,\xi ,h),h)\equiv G(\rrr ,\xi ,h) \).

\item If \( h^{\prime }\leq _{T}g \) then
\( G(\mu ,G(\nu ,\xi ,h),g)\equiv G(\rrr ,G(\mu ,\xi ,g),h) \).

\item If \( g^{\prime }\leq _{T}h \) then
\( G(\mu ,G(\nu ,\xi ,h),g)\leq G(\mu \oplus \nu ,\xi ,h) \).
 \end{enumerate}
\end{Proposition}

The next assertion, proved by arguments close to those in \cite{s83},  is formally new here.

\begin{Proposition}\label{p-r2}
\begin{enumerate}
\item If \( \mu \leq ^{h^{\prime }}\nu  \) and \( \nu  \) is \( h^{\prime } \)-complete
then \( G(\mu \oplus \mu _{1},\nu ,h)\equiv G(\mu _{1},\nu ,h) \).

\item If \( \nu  \) is \( h^{\prime } \)-complete then
\( G(\rrr \oplus \mu _{1},\nu ,h)\equiv G(\mu \oplus \mu _{1},\nu ,h) \).

\item If  \( \mu  \) (resp. \( \nu  \)) is \( h^{\prime } \)-complete
then \( \lambda \nu .\rrr  \) (resp. \( \lambda \mu .\rrr  \)) is a closure operator
on  \( (^{\omega }S;\leq ) \)  (resp. on \( (^{\omega }S;\leq ^{h}) \)).

\item If \( G(\mu ,\nu ,h)\leq^h G(\mu_1 ,\nu_1 ,h) \) and $\nu$ is $h'$-precomplete then $\nu\leq\nu_1$ or \( G(\mu ,\nu ,h)\leq\mu_1 \).

 \end{enumerate}
\end{Proposition}

{\em Proof.} (1) By Proposition {\ref{p-r}(3)}, it  suffices to show that \( G(\mu \oplus
\mu _{1},\nu ,h)\leq G(\mu _{1},\nu ,h) \).
Let \( u \) be an \( h^{\prime } \)-computable function with \( \mu =\nu \circ u \),
\( \omega _{i}=\{2x+i\mid x<\omega \} \) for $i<2$, and let
\[
f_{1}\langle m,x\rangle =\left\{ \begin{aligned}
m,&\; \text {if}\; \tilde{\varkappa} ^{h}(x)\uparrow \vee
\tilde{\varkappa} ^{h}(x)\downarrow \in \omega _{1},\\
u(\frac{\tilde{\varkappa} ^{h}(x)}{2}),&\; \text {if}\;
\tilde{\varkappa} ^{h}(x)\downarrow \in \omega _{0}.
\end{aligned}\right.
\]
Then \( f_{1}\leq _{T}h^{\prime } \), hence \( \nu \circ f_{1}=\nu \circ f \)
for some computable function \( f \).
Let \( g \) be a computable function such that
\[
\tilde{\varkappa} ^{h}g(x)=\left\{ \begin{aligned}
\uparrow, & \, \text {if}\, \tilde{\varkappa} ^{h}(x)\uparrow
\vee \tilde{\varkappa} ^{h}(x)\downarrow \in \omega _{0},\\
u(\frac{\tilde{\varkappa} ^{h}(x)-1}{2}), &\, \text {if}\,
\tilde{\varkappa} ^{h}(x)\downarrow \in \omega _{1}.
\end{aligned}\right.
\]
Then the computable function \( \langle m,x\rangle \mapsto
\langle f(x),g(x)\rangle  \) reduces \( G(\mu \oplus \mu _{1},\nu ,h) \)
to \( G(\mu _{1}\nu ,h) \).

(2) By Proposition {\ref{p-r}(3)}, it suffices to check that \( G(\rrr \oplus
\mu _{1},\nu ,h)\leq G(\mu \oplus \mu _{1},\nu ,h) \).
Let
\[
f_{1}\langle m,x\rangle =\left\{ \begin{aligned}
m, &\, \text {if}\, \tilde{\varkappa} ^{h}(x)\uparrow \vee
\tilde{\varkappa} ^{h}(x)\downarrow  \in \omega _{1},\\
\frac{\tilde{\varkappa} h(x)}{2}, &\, \text {if}\,
\tilde{\varkappa} ^{h}(x)\downarrow \in \omega _{0}.
\end{aligned}\right.
\]
Then \( f_{1}\leq _{T}h^{\prime } \), hence \( \nu \circ f_{1}=\nu \circ f \)
for some computable function \( f \). Let \( g \) be a computable function satisfying
\[
\tilde{\varkappa} ^{h}g(x)=\left\{ \begin{aligned}
\uparrow,& \, \text {if}\, \tilde{\varkappa} ^{h}(x)\uparrow \vee (\tilde{\varkappa} ^{h}(x)\downarrow
\in \omega _{0}\wedge \tilde{\varkappa} ^{h}r(\frac{\tilde{\varkappa} ^{h}(x)}{2})\uparrow ),\\
2\tilde{\varkappa} ^{h}r(\frac{\tilde{\varkappa} ^{h}(x)}{2}),&\, \text {if}\, (\tilde{\varkappa} ^{h}(x)\downarrow
\in \omega _{0}\wedge \tilde{\varkappa} ^{h}r(\frac{\tilde{\varkappa} ^{h}(x)}{2})\downarrow ),\\
\tilde{\varkappa} ^{h}(x),&\, \text {if}\, \tilde{\varkappa} ^{h}(x)\downarrow \in \omega _{1}.
\end{aligned}\right.
\]
Then the computable function \( \langle m,x\rangle \mapsto \langle f\langle m,x\rangle ,
g(x)\rangle  \) reduces \( G(\rrr \oplus \mu _{1},\nu ,h) \)
to \( G(\mu \oplus \mu _{1},\nu ,h) \).

(3) Recall that a closure operator on a preorder $(Q;\leq)$ is a monotone function $f:Q\to Q$ such that $x\leq f(x)$ and $f(f(x))\leq f(x)$.
The monotonicity and the  property $x\leq f(x)$ for the operators in formulation
are contained in Proposition {\ref{p-r}}. It remains
to check that \( G(\mu ,\rrr ,h)\leq \rrr  \) for an
 \( h^{\prime } \)-complete \( \mu  \),
and \( G(\rrr ,\nu ,h)\leq^h \rrr
\) for an
 \( h^{\prime } \)-complete \( \nu  \). By Proposition {\ref{p-r}}(2),
 \( G(\mu ,\mu ,h)\leq ^{h^{\prime }}\mu  \).
Since \( \mu  \) is \( h^{\prime } \)-complete, we get \( G(\mu ,\mu ,h),\nu ,h)\leq \mu  \).
By item (1) and Proposition  {\ref{p-r}}(3) we obtain
\( G(\mu ,\rrr ,h)\equiv G(G(\mu ,\mu ,h),\nu ,h)\leq \rrr.
\)

It remains to check the assertion for \( \nu  \). As above,
 \( G(\nu ,\nu ,h)\leq \nu  \).
By  (1) and Proposition  {\ref{p-r}}(5),
\( G(\rrr ,\nu ,h)\equiv G(\mu ,G(\nu ,\nu ,h),h)\leq \rrr
\).
\medskip

(4) Let $f$ be an $h$-computable function that reduces \( G(\mu ,\nu ,h)\) to \(G(\mu_1 ,\nu_1 ,h) \) and let
 $$Y_x=\{y\mid \tilde{\varkappa}^hrf\langle \varkappa^h_{l(x)}(y),\varkappa^h_{r(x)}(y)\rangle\downarrow\}.$$
  Since the set $Y_x$ is $h$-c.e. uniformly on $x$, there are uniform sequences $\{\varphi_x\},\{\psi_x\}$ of $h$-c.p. functions such that $\psi_x$ is a bijection between a (unique) initial segment of $\omega$ and $Y_x$, and $\varphi_x=\psi_x^{-1}$. Let $\bar{\varphi}_x(y)=\varphi_x(y)$ for $y\in Y_x$ and $\bar{\varphi}_x(y)=y$ otherwise, then $\bar{\varphi}_x\leq_Th'$ uniformly on $x$. Since $\nu$ is $h'$-precomplete, there is a uniform sequence $\{u_x\}$ of computable total functions such that $\nu\circ l\circ\bar{\varphi}_x=\nu\circ u_x$. Since $\tilde{\varkappa}^h$ is $h$-complete, there is a uniform sequence $\{v_x\}$ of total computable functions such that $\tilde{\varkappa}^h\circ r\circ\varphi_x=\tilde{\varkappa}^h\circ v_x$. By the double recursion theorem, there is $e$ such that $u_e=\varkappa_{l(e)}$ and $v_e=\varkappa_{r(e)}$.

It suffices to show that if $Y_e$ is finite then $\nu\leq\nu_1$, else \( G(\mu ,\nu ,h)\leq\mu_1 \). Let first $Y_e$ be finite. For any $z=\langle z_1,z_2\rangle\in\omega\setminus Y_e$ with $\tilde{\varkappa}^h(z_2)\uparrow$ we have $\nu(z_1)=G(\mu ,\nu ,h)(z)$. Since $\varphi_e(z)\uparrow$, we have $\tilde{\varkappa}^hv_e(z)\uparrow$, hence $G(\mu ,\nu ,h)\langle u_e(z),v_e(z)\rangle=\nu u_e(z)=\nu l\bar{\varphi}_e(z)=\nu l(z)=\nu(z_1)$. Therefore,
 $$\nu(z_1)=G(\mu ,\nu ,h)\langle u_e(z),v_e(z)\rangle=G(\mu ,\nu ,h)\langle \varkappa_{l(e)}(z),\varkappa_{r(e)}(z)\rangle=\nu_1 lf\langle u_e(z),v_e(z)\rangle,
 $$
 so $\nu\leq\nu_1$.
 
Let now $Y_e$ be infinite, then $\varphi_e:Y_e\to\omega$ is a bijection and $\varphi_e\psi(m)=m$ for every $m<\omega$. Setting $y=\psi_e(m)$, we obtain $\nu l\varphi_e(y)=\nu u_e(y)$ and $\mu\tilde{\varkappa}^h\varphi_e(y)=\mu\tilde{\varkappa}^hv_e(y)$, hence
 \begin{eqnarray*}
 G(\mu ,\nu ,h)(m)=G(\mu ,\nu ,h)\varphi_e(y)=G(\mu ,\nu ,h)\langle u_e(y),v_e(y)\rangle=\\ G(\mu_1 ,\nu_1 ,h)\langle f \varkappa_{l(e)}(y),\varkappa_{r(e)}(y)\rangle=\mu_1\tilde{\varkappa}^hf\langle u_e\psi_e(m),v_e\psi_e(m)\rangle.
 \end{eqnarray*}
 Thus, \( G(\mu ,\nu ,h)\leq\mu_1 \).
 \qed

\subsection{An algebra of $k$-partitions of $\omega$}\label{algebra}

For $S=\bar{k}$, $2\leq k\leq\omega$, the operation $G$ from Definition \ref{operG} becomes  a ternary operation on $k^\omega$. Let  $A_k$ be the subalgebra of $(k^\omega;\oplus,G)$ generated by the constant functions $\bi=\lambda x.i$, $i<k$.  Here we prove a series of facts about this subalgebra which relates it to the structure of iterated $k$-forests. 

Let $\mathbb{T}_k$ be the set of variable-free terms of signature $\sigma_k=\{i,\oplus,G\mid i<k\}$, then $A_k=\{\bu\mid u\in \mathbb{T}_k\}$ where $\bu$ is the value of $u$ in $(k^\omega;\sigma_k)$.
Let $\bo^{(0)}=\bo$ and $\bo^{(n+1)}=G(\bo,\boo,\bo^{(n)})$, then $\{\bo^{(n)}\}$ essentially coincides with the usual  iterations of Turing jump starting with $\bo$, see  Proposition \ref{p-r}(7). The following fact characterizes the quotient-poset of $(A_k;\leq_T)$. 
\medskip

{\bf Fact 1.} Given $u\in \mathbb{T}_k$, one can compute $n=n(u)$ with $\bu\equiv_T\bo^{(n)}$. Therefore, $A_k\equiv_T\{\bo^{(n)}\mid n<\omega\}$ and hence $(A_k;\leq_T)\simeq(\omega;\leq)$.
\medskip

{\em Proof.} Define $n(u)$ by induction on (the rank of) $u$ as follows: $n(i)=0$ for every $i<k$, $n(u_1\oplus u_2)=max(n(u_1),n(u_2))$, $n(G(u_1,u_2,u_3))=0$ if $u_1=u_2=i$ for some $i<k$ and $n(G(u_1,u_2,u_3))=max(n(u_1),n(u_2),n(u_3)+1)$ otherwise. By  Proposition \ref{p-r}(7), the function $u\mapsto n(u)$ works.
 \qed
 
 For any $n<\omega$ we define the binary operation $\cdot^n$ on $k^\omega$ by $\nu\cdot^n\mu=G(\mu,\nu,\bo^{(n)})$; note that $A_k$ is closed under these operations. Let $\mathbb{T}^*_k$ be the set of variable-free terms of signature $\sigma^*_k=\{i,\oplus,\cdot^n\mid  i<k, n<\omega\}$. 
\medskip

{\bf Fact 2.}  Given $u\in \mathbb{T}_k$, one can compute $u^*\in \mathbb{T}^*_k$ with $\bu\equiv\mathbf{u^*}$, and vice versa. In particular,  $A_k\equiv\{\bu\mid u\in \mathbb{T}^*_k\}$.
\medskip

{\em Proof.} Define $u\mapsto u^*$ by induction on  $u$ as follows: $i^*=i$ for every $i<k$, $(u_1\oplus u_2)^*=u_1^*\oplus u_2^*$, and $G(u_1,u_2,u_3)^*=u_2^*\cdot^nu_1^*$ where $n=n(u_3)$ is from Fact 1. By items (3,5,6) of Proposition \ref{p-r}, $\bu\equiv\mathbf{u^*}$ for every $u\in \mathbb{T}_k$. The opposite direction (computing from any given $u\in \mathbb{T}^*_k$ some $v\in \mathbb{T}_k$ with $\bu\equiv\mathbf{v}$) is considered similarly.
 \qed
 \medskip
 
From now on we will work with the signature $\sigma^*_k$ and its subsignatures. We abbreviate $\leq^{\bo^{(n)}}$ to $\leq^n$, so in particular  $\leq^0=\leq$. 
For any $n<\omega$, let $A^n_k=\{\nu\in A_k\mid\nu\leq_T\bo^{(n)}\}$. By Fact 1, $A^n_k=\{\bu\mid n(u)\leq n\}$. We give a more constructive characterization of $A^n_k$. Let $\mathbb{T}^n_k$ be the set of variable-free terms of  $\sigma^n_k=\{i,\oplus,\cdot^p\mid i<k, p<n\}$, then $\bigcup_n\mathbb{T}^n_k=\mathbb{T}^*_k$.

 \medskip
{\bf Fact 3.}  We have: $A^n_k\equiv\{\bu\mid u\in \mathbb{T}^n_k\}$.
\medskip

{\em Proof.} The inclusion $\supseteq$ follows from the proof of Fact 1. For the converse, it suffices to prove by induction on $u\in \mathbb{T}^*_k$ that if $\bu\in A^n_k$ then $\bu\equiv\bu^*$ for some $u^*\in \mathbb{T}^n_k$. If $u=i$ or $u=u_1\oplus u_2$, we respectively set $u^*=i$ or $u^*=u^*_1\oplus u^*_2$. Now let $u=u_1\cdot^pu_2$. If $p<n$, we can set $u^*=u^*_1\cdot^p u^*_2$. Finally, let $n\leq p$. If $u_1=u_2=i$ for some $i<k$, we set $u^*=i$. The remaining case is not possible because in this case $n(u)\geq p+1>n$ by Fact 1, hence $\bo^{(n+1)}\leq_T\bu$ and $\bu\not\in A^n_k$.
 \qed

We define the family $\{f^n_m\}_{n<\omega}$ of functions $f^n_m:\mathcal{T}_m(\bar{k})\to \mathbb{T}^*_k$ by induction on $m$ as follows. Let $f^n_0(i)=i$ for all $i<k,n<\omega$. It remains to define $f^n_{m+1}$ from $f^{n+1}_m$. Let $(T,t)\in\mathcal{T}_k(m+1)=\mathcal{T}_{\mathcal{T}_k(m)}$. If $T$ is singleton we set $f^n_{m+1}(T)=f^{n+1}_m(t(\varepsilon))$, otherwise  we set
 $$f^n_{m+1}(T)=f^{n+1}_m(t(\varepsilon))\cdot^n(\bigoplus\{f^n_{m+1}(T(i))\mid i\in\omega\cap T\}),$$
  using induction on the rank of $T$. We also define  functions $\bff^n_m:\mathcal{T}_m(\bar{k})\to A_k$ by $\bff^n_m=ev\circ f^n_m$ where $ev:\mathbb{T}^*_k\to A_k$ is the evaluation function $ev(u)=\bu$.
\medskip

{\bf Fact 4.}  For  $T,V\in\mathcal{T}_k(m)$ we have: $\bff^n_m(T)$ is $\bo^{(n)}$-complete, and $T\leq_hV$ iff $\bff^n_m(T)\leq \bff^n_m(V)$.
\medskip

{\em Proof.}  The $\bo^{(n)}$-completeness of $\bff^n_m(T)$ follows by induction on $m$ from the definition and Proposition \ref{p-r}(9). The second fact is checked by induction on $m$ and the ranks of trees $(T,t)$ and $(V,v)$. For $m=0$ the assertion is obvious, so let $T,V\in\mathcal{T}_k(m+1)$. Let first  $T$ be singleton, then $\bff^n_{m+1}(T)=\bff^{n+1}_m(t(\varepsilon))$, which we temporarily denote by $\nu$, is $\bo^{(n+1)}$-complete. By the definition, $\bff^n_{m+1}(V)$ is the value of some $\{\oplus,\cdot^n\}$-term $u(x_1,\ldots,x_p)$ whose variables take values in $\{f^{n+1}_m(t(\sigma))\mid \sigma\in V\}$. By Propositions \ref{hprecomp}(3) and \ref{p-r2}(4), if $\nu\leq^n\nu_1\cdot^n\mu_1$ ($\nu\leq^n\nu_1\oplus\mu_1$) then $\nu\leq\nu_1$ or $\nu\leq\mu_1$. Therefore, $T\leq_hV$ iff $t(\varepsilon)\leq_hv(\sigma)$ for some $\sigma\in V$ iff $\bff^{n+1}_m(t(\varepsilon))\leq^n\bff^{n+1}_m(v(\sigma))$ for some $\sigma\in V$ iff $\bff^n_{m+1}(T)\leq \bff^n_{m+1}(V)$. 

Now let $V$ be singleton, then $\bff^n_{m+1}(V)=\bff^{n+1}_m(v(\varepsilon))$, which we temporarily denote by $\mu$, is $\bo^{(n+1)}$-complete. By the definition, $\bff^n_{m+1}(T)$ is the value of some $\{\oplus,\cdot^n\}$-term $u(x_1,\ldots,x_p)$ whose variables take values in $\{\bff^{n+1}_m(t(\tau))\mid \tau\in T\}$. By Proposition \ref{p-r}(8) we have $\mu\cdot^n\mu\equiv\mu$ (and of course also $\mu\oplus\mu\equiv\mu$). Therefore, $T\leq_hV$ iff $t(\tau)\leq_hv(\varepsilon)$ for all $\tau\in T$ iff $\bff^{n+1}_m(t(\tau))\leq\mu$ for all $\tau\in T$ iff $\bff^n_{m+1}(T)\leq \bff^n_{m+1}(V)$.

Finally, let both $T$ and $V$ be non-singletons. If $t(\varepsilon)\leq v(\varepsilon)$ then $\bff^{n+1}_m(t(\varepsilon))\leq\bff^{n+1}_m(v(\varepsilon))$ and therefore we have: $T\leq_hV$ iff  $T(i)\leq V$ for all $i\in\omega\cap T$ iff $\bff^n_{m+1}(T(i))\leq \bff^n_{m+1}(V)$ for all $i\in\omega\cap T$ iff $\bff^n_{m+1}(T)\leq \bff^{n+1}_m(v(\varepsilon))\cdot^n  \bff^n_{m+1}(V)$  iff $\bff^n_{m+1}(T)\leq \bff^n_{m+1}(V)$ (the latter equivalence uses Proposition \ref{p-r}(10)).
If $t(\varepsilon)\not\leq v(\varepsilon)$ then $\bff^{n+1}_m(t(\varepsilon))\not\leq \bff^{n+1}_m(v(\varepsilon))$ and therefore we have: $T\leq_hV$ iff  $T\leq V(j)$ for some $j\in\omega\cap V$ iff $\bff^n_{m+1}(T)\leq \bff^n_{m+1}(V(j))$ for some $j\in\omega\cap V$  iff $\bff^n_{m+1}(T)\leq \bff^n_{m+1}(V)$ (the latter equivalence uses Proposition \ref{p-r2}(4)).
 \qed

Using Proposition \ref{slattice}, we extend the function $f^n_m:\mathcal{T}_m(\bar{k})\to \mathbb{T}^*_k$  to a function $f^n_m:\mathcal{F}_m(\bar{k})\to \mathbb{T}^*_k$ (denoted for simplicity by the same name) by $f^n_m(T_0\sqcup\cdots\sqcup T_p)=f^n_m(T_0)\oplus\cdots\oplus f^n_m(T_p)$ where $T_i$ are trees.  We also define  functions $\bff^n_m:\mathcal{F}_m(\bar{k})\to A_k$ by $\bff^n_m=ev\circ f^n_m$, as above. In the next fact we use the operation $\cdot$ on forests defined at the end of Section \ref{trees}.

 \medskip
 
{\bf Fact 5.}  For all $F,G\in\mathcal{F}_m(\bar{k})$ we have: $F\leq_hG$ iff $\bff^n_m(F)\leq^n \bff^n_m(G)$, $\bff^n_m(F\sqcup G)\equiv^n\bff^n_m(F)\oplus \bff^n_m( G)$, and $\bff^n_m(F\cdot G)\equiv^n\bff^n_m(F)\cdot^n \bff^n_m( G)$.
\medskip

{\em Proof.} The first and second assertions follow from Proposition \ref{slattice} and Fact 4.
For the third assertion, we use induction on the  cardinality of $F$. Let first $F=T$ be a tree. If $T$ is singleton, the assertion holds by the definition of $\bff^n_m$. Otherwise, using the induction, the associativity of $\cdot$ (Proposition \ref{iterh}(2)) and of $\cdot^n$ (which holds by Proposition \ref{p-r}(10)) and abbreviating $\bff^n_m$ to $\bff$, we obtain:
 \begin{eqnarray*}
 \bff(T\cdot G)=\bff((\{\varepsilon\}\cdot(T\setminus\{\varepsilon\})\cdot G)=\bff(\{\varepsilon\}\cdot(T\setminus\{\varepsilon\}\cdot G))=\bff(\{\varepsilon\})\cdot^n\bff(T\setminus\{\varepsilon\})\cdot G)\equiv^n\\ \bff(\{\varepsilon\})\cdot^n(\bff(T\setminus\{\varepsilon\})\cdot^n\bff( G))\equiv^n(\bff(\{\varepsilon\})\cdot^n\bff(T\setminus\{\varepsilon\})\cdot^n\bff( G)\equiv^n\bff(T)\cdot^n\bff(G).
 \end{eqnarray*}

Let now $F$ be not a tree, then $F\equiv_hF_1\sqcup F_2$ for some $F_1,F_2\in\mathcal{F}_m(\bar{k})$ of lesser cardinalities than $F$. Using the induction, Facts 5 and 7,  the right distributivity of $\cdot$ w.r.t. $\sqcup$ (see the end of Section \ref{trees}) and of $\cdot^n$ w.r.t. $\oplus$ (which is essentially Proposition \ref{p-r}(1)), we obtain:
 \begin{eqnarray*}
 \bff(F\cdot G)\equiv^n\bff((F_1\sqcup F_2)\cdot G)\equiv^n\bff((F_1\cdot G)\sqcup(F_2\cdot G))\equiv^n\bff(F_1\cdot G)\oplus \bff(F_2\cdot G)\equiv^n\\(\bff(F_1)\cdot^n\bff( G))\oplus (\bff(F_2)\cdot^n\bff( G))\equiv^n(\bff(F_1)\oplus \bff(F_2))\cdot^n\bff( G)\equiv^n\bff(F)\cdot^n\bff( G).
 \end{eqnarray*} 
 \qed

We define the unary functions $s$ and $r$ on $\mathbb{T}^*_k$ by induction on terms as follows: $s(i)=r(i)=i$, $s(u_1\oplus u_2)=s(u_1)\oplus s(u_2)$, $r(u_1\oplus u_2)=r(u_1)\oplus r(u_2)$,  $s(u_1\cdot^p u_2)=s(u_1)\cdot^{p+1} s(u_2)$,  $r(u_1\cdot^0 u_2)=r(u_1)\oplus r(u_2)$, $r(u_1\cdot^{p+1} u_2)=r(u_1)\cdot^{p} r(u_2)$. The next fact is obvious (the last assertion holds because the term $f^n_m(F)$ has no entries of $\cdot^0$). 

\medskip

{\bf Fact 6.} We have: $r(s(u))=u$, $s(\mathbb{T}_k^n)\subseteq\mathbb{T}_k^{n+1}$,   $r(\mathbb{T}_k^{n+1})\subseteq\mathbb{T}_k^n$, and $s(r(f^n_m(F)))=f^n_m(F)$ for $n>0$.

The next two facts describe relationships of  operations $s$ and $r$ (on the labeled forests and on the terms) to the functions $f^n_m$. 

\medskip
{\bf Fact 7.} For any $F\in\mathcal{F}_m(\bar{k})$ we have: $s(f^n_m(F))=f^{n+1}_m(F)$ and  $r(f^{n+1}_m(F))=f^{n}_m(F)$.
\medskip

{\em Proof.} Since all $s,r,f^n_m$ respect $\sqcup$ and $\oplus$, it suffices to check this for the case when $F=(T,t)$ is a tree. We argue by induction on $m$. For $m=0$ the assertion is clear since $F=i<k$, so let $m>0$.  If $T$ is singleton then $s(f^n_m(F))=s(f^{n+1}_{m-1}(t(\varepsilon)))=f^{n+2}_{m-1}(t(\varepsilon))=f^{n+1}_{m}(F)$ and  $r(f^{n+1}_m(F))=r(f^{n+2}_{m-1}(t(\varepsilon)))=f^{n+1}_{m-1}(t(\varepsilon))=f^{n}_m(F)$. Otherwise, by induction on the rank of $T$ we have:
 $$s(f^n_m(F))=s(f^{n}_{m}(t(\varepsilon))\cdot^n\bigoplus_if^n_m(T(i)))=f^{n+1}_{m}(t(\varepsilon))\cdot^{n+1}\bigoplus_if^{n+1}_m(T(i))=f^{n+1}_{m}(F)$$
  and  $r(f^{n+1}_m(F))=r(f^{n+1}_{m}(t(\varepsilon))\cdot^n\bigoplus_if^{n+1}_m(T(i)))$. This equals to $f^{n}_{m}(t(\varepsilon))\oplus\bigoplus_if^{n}_m(T(i))$ for $n=0$ and to $f^{n}_{m}(t(\varepsilon))\cdot^{n-1}\bigoplus_if^{n}_m(T(i))$ for $n>0$; in any case, it  equals to $f^{n}_m(F)$.
 \qed
 
 \medskip
  
{\bf Fact 8.}  For all $F\in\mathcal{F}_m(\bar{k})$ and $G\in\mathcal{F}_{m+1}(\bar{k})$ we have: $s( f^n_m(F))=f^n_{m+1}(s(F))$, $r(f^0_{m+1}(G))=f^0_{m}(r(G))$, and $ev(r(f^n_{m+1}(G)))\equiv ev(f^n_{m}(r(G)))$ for $n>0$.
\medskip

{\em Proof.} Since all involved functions respect $\sqcup$ and $\oplus$, for the first equality we have to show that $s(f^n_m(F))=f^n_{m+1}(s(F))$ for every tree $F=(T,t)\in\mathcal{T}_m(\bar{k})$, for the second equality we have to show that $r(f^0_{m+1}(F))=f^0_{m}(r(F))$ for every tree $G=(T,t)\in\mathcal{T}_{m+1}(\bar{k})$, and similarly for the third assertion. 
We argue by induction on $m$. For $m=0$ the first equality is clear since $F=i<k$, so let $m>0$. 
Since $f^n_{m+1}(s(F))=f^{n+1}_m(F)$, for the first equality it suffices to show that $s(f^n_m(F))=f^{n+1}_m(F)$; this holds by Fact 7. 

The second and third equalities are checked by induction on the rank of $T$. If $T$ is singleton, i.e. $T=s(V)$ for some $V\in\mathcal{T}_m(\bar{k})$ then $r(f^n_{m+1}(T))=r(f^{n+1}_{m}(V))$ and $f^n_m(r(T))=f^n_m(r(s(V)))=f^n_m(V)$, so the equalities hold by Fact 7. If $T$ is not singleton then we have
\begin{eqnarray*} 
 f^n_{m}(r(T))=f^n_m(r(t(\varepsilon)\sqcup\bigsqcup_iT_i))=f^n_{m}(r(t(\varepsilon)))\oplus\bigoplus_if^n_{m}(r(T_i))\\
  \text{and } r(f^n_{m+1}(T))=r(f^{n+1}_{m}(t(\varepsilon)))\cdot^n\bigoplus_if^n_{m+1}(T(i)).
\end{eqnarray*}  
   The latter expression equals to 
  $r(f^n_{m+1}(t(\varepsilon)))\oplus\bigoplus_ir(f^n_{m+1}(T(i)))$ 
 (hence to $f^n_{m}(r(T))$, yielding the second equality) for $n=0$ and to $r(f^n_{m+1}(t(\varepsilon)))\cdot^{n-1}\bigoplus_ir(f^n_{m+1}(T(i)))$ for $n>0$. For $n>0$, by Proposition \ref{p-r}(2) we get
 $$ev(r(f^n_{m+1}(t(\varepsilon)))\cdot^{n-1}\bigoplus_ir(f^n_{m+1}(T(i)))\equiv^nev(r(f^n_{m+1}(t(\varepsilon)))\oplus\bigoplus_ir(f^n_{m+1}(T(i))))$$
 which equals to the $\bo^n$-complete numbering $ev(f^n_m(r(T)))$; this yields the third equality.
 \qed

The assignment $\bff^n_m(F)\mapsto\bff^n_{m+1}(s(F))$ defines a multifunction from $B^n_m=\bff^n_m(\mathcal{F}_m(\bar{k}))$ into $B^n_{m+1}$ which is denoted by $s$. From Fact 8 it follows that $s$ induces  a function on the $\equiv^n$-classes. Similarly, let $r:B^n_{m+1}\to B^n_{m}$ be the function on $\equiv^n$-classes corresponding to the assignment  $\bff^n_{m+1}(F)\mapsto\bff^n_{m}(r(F))$.

 \medskip
 
{\bf Fact 9.}  For all  $F,G\in\mathcal{F}_m(\bar{k})$ and $p<m$ we have: 

$\bff^n_m(F)\leq^n\bff^n_m(G)$ iff $\bff^n_{m+1}(F)\leq^{n+1}\bff^n_{m+1}(G)$ iff $\bff^n_{m+1}(F)\leq^0\bff^n_{m+1}(G)$; 

$\bff^n_{m}(F)\cdot^{n+p+1}\bff^n_{m}(G)\equiv^n\bff^n_{m}(sr(F))\cdot^{n+p+1}\bff^n_{m}(sr(G))$; 

$\bff^n_m(F)\leq^{n+p+1}\bff^n_m(G)$ iff $\bff^n_{m-1}(r(F))\leq^{n+p}\bff^n_{m-1}(r(G))$.
\medskip

{\em Proof.} Let $F=T_0\sqcup\cdots\sqcup T_q$ and $G=V_0\sqcup\cdots\sqcup V_l$ be minimal forests  decomposed to trees. Then $F\leq_hG$ iff $s(F)\leq_h's(G)$ iff $s(F)\leq_hs(G)$ by Section \ref{trees}, and $\bff^n_{m+1}(F)=\bigoplus_i\bff^{n+1}_{m}(T_i)$ and $\bff^n_{m+1}(G)=\bigoplus_j\bff^{n+1}_{m}(V_j)$ by the definition of $\bff^n_{m+1}(F)$. Since any $\bff^{n+1}_{m}(V_j)$ is $\mathbf{0}^{(n+1)}$-complete by Fact 4, the first assertion follows from Proposition \ref{hprecomp}.

 For $n>0$ the second assertion follows from Fact 6 (because $s(r(f^n_m(F)))=f^n_m(F)$ and similarly for $G$), so let $n=0$. Let $t$ and $g$ be the labelings of $F$ and $G$, resp. By the definition of $f^0_m$, $f^0_m(F)$ (resp. $f^0_m(G)$) is a term of signature $\{\oplus,\cdot^0\}$ from $f^1_{m-1}(t(\tau))$, $\tau\in F$ (resp. from $f^1_{m-1}(g(\sigma))$, $\sigma\in G$). By items (2) and (12) of Proposition \ref{p-r},
  $$\bff^0_{m}(F)\cdot^{p+1}\bff^0_{m}(G)\equiv^0(\bigoplus_\tau\bff^1_{m-1}(t(\tau)))\cdot^{p+1}\bigoplus_\sigma\bff^1_{m-1}(g(\sigma))).$$
   On the other hand, since $sr(f^1_{m-1}(t(\tau)))=f^1_{m-1}(t(\tau))$, $srf^0_{m}(F)$ (resp. $srf^0_{m}(G)$) is a term of signature $\{\oplus\}$ from $f^1_{m-1}(t(\tau))$, $\tau\in F$ (resp. from $f^1_{m-1}(g(\sigma))$, $\sigma\in G$). Therefore,
    $$sr(\bff^0_{m}(F))\cdot^{p+1}sr(\bff^0_{m}(G))\equiv^0(\bigoplus_\tau\bff^1_{m-1}(t(\tau)))\cdot^{p+1}\bigoplus_\sigma\bff^1_{m-1}(g(\sigma))),$$
     completing the proof of the second assertion.
 
By the first assertion, the third assertion follows from: $\bff^n_m(F)\leq^{n+p+1}\bff^n_m(G)$ iff $\bff^n_m(sr(F))\leq^{n+p+1}\bff^n_m(sr(G))$. This follows from the argument of the preceding paragraph (for $n$ in place of $0$) which shows that $\bff^n_m(F)\equiv^{n+p+1}\bff^n_m(sr(F))$, and similarly for $G$.
 \qed

 \medskip
 
{\bf Fact 10.}  For all  $F,G\in\mathcal{F}_m(\bar{k})$ and $p<m$ we have: $\bff^n_m(F\cdot^p G)\equiv^n\bff^n_m(F)\cdot^{n+p} \bff^n_m( G)$.
\medskip

{\em Proof.} We argue by induction on $p$. For $p=0$ the assertion holds by Fact 5, so let $p+1<m$ and jet the assertion hold for $\cdot^p$. Using the functions $s,r$ from Section \ref{trees}, the corresponding functions on $\equiv^n$-classes defined above, Facts 8 and 9, the induction hypothesis, and Proposition \ref{iterh}(2), we obtain:
 \begin{eqnarray*}
  \bff^n_m(F\cdot^{n+p+1} G)=\bff^n_m(s(r(F)\cdot^{n+p}r(G)))\equiv^ns(\bff^n_{m-1}(r(F)\cdot^{n+p}r(G)))
  \equiv^n\\s(\bff^n_{m-1}(r(F))\cdot^{n+p}\bff^n_{m-1}(r(G))))\equiv^ns(r(\bff^n_{m}(F))\cdot^{n+p}r(\bff^n_{m}(G)))\equiv^n\\ s(r(\bff^n_{m}(F)))\cdot^{n+p+1}s(r(\bff^n_{m}(G))))\equiv^n\bff^n_{m}(F)\cdot^{n+p+1}\bff^n_{m}(G).
 \end{eqnarray*}
 \qed

 \medskip
 
{\bf Fact 11.}  The function  $\bff^n_m$ induces an isomorphism between the quotient-structure of the structure  $(\mathcal{F}_m(\bar{k});\sqcup,\cdot^0,\ldots,\cdot^{m-1},\leq^0,\ldots,\leq^m)$ under $\equiv_h$ and  the quotient-structure of $(A^{n+m}_k;\oplus,\cdot^n,\ldots,\cdot^{n+m-1},\leq^n,\ldots,\leq^{n+m})$ under $\equiv^n$.
\medskip

{\em Proof.} By Fact 5, $\bff^n_m:\mathcal{F}_m(\bar{k})\to B^n_m$ is a semilattice embedding of $(\mathcal{F}_m(\bar{k});\sqcup,\leq_h)$ into $(A_k;\oplus,\leq^n)$. By the definition of $f^n_m$, any $f^n_m(F)$ is a term of signature $\{\oplus,\cdot^0,\ldots,\cdot^{m-1}\}$ from $i<k$ (in case $m=0$ the multiplications are absent). By the definition of $\bff^n_m$, any $\bff^n_m(F)$ is a the value of a variable-free term of signature $\{i,\oplus,\cdot^n,\ldots,\cdot^{n+m-1}\mid i<k\}$. By Fact 3, $B^n_m\subseteq A^{n+m}_k$, i.e. $\bff^n_m:\mathcal{F}_m(\bar{k})\to A^{n+m}_k$.

We show that $B^n_m\equiv^n A^{n+m}_k$, i.e. $\bff^n_m$ is an isomorphism between the quotient-structures of $(\mathcal{F}_m(\bar{k});\sqcup,\leq_h)$ and $(A^{n+m}_k;\oplus,\leq^n)$. Let $\mathbb{T}^{n,m}_k$ be the set of variable-free terms of signature $\sigma^{n,m}_k=\{i,\oplus,\cdot^p\mid i<k, n\leq p<n+m\}$. Note that $\mathbb{T}^{0,m}_k=\mathbb{T}^{m}_k$ and $\mathbb{T}^{n,m}_k\subseteq\mathbb{T}^{n+m}_k$. By Fact 3, it suffices to show that, given $u\in \mathbb{T}^{n+m}_k$, one can compute $u^*\in \mathbb{T}^{n,m}_k$ with $\bu\equiv^n\mathbf{u^*}$. We define $u\mapsto u^*$ by induction on  $u$ as follows: $i^*=i$ for every $i<k$, $(u_1\oplus u_2)^*=u_1^*\oplus u_2^*$,  $(u_1\cdot^pu_2)^*=u_1^*\cdot^pu_2^*$ for $p\geq n$, and $(u_1\cdot^pu_2)^*=u_1^*\oplus u_2^*$ for $p< n$. For the first two cases the assertion is clear, so it suffices to check that $(\bu_1\cdot^p\bu_2)^*\equiv^n\bu_1^*\cdot^p\bu_2^*$. By induction we have $\bu_1\equiv^n\bu_1^*$ and $\bu_2\equiv^n\bu_2^*$. For $p\geq n$ we get $\bu_1^*\cdot^p\bu_2^*\equiv^n(\bu_1\cdot\bu_2)^*$ by items (3) and (5) of Proposition \ref{p-r}. For $p<n$ we have $\bu_1^*\cdot^p\bu_2^*\equiv^n\bu_1^*\oplus\bu_2^*=$ by Proposition \ref{p-r}(2) which imply the desired property.

By Fact 10, the  function $\bff^n_m$ respects the corresponding multiplications, so it remains to show that it also respects the corresponding preorders. This is by induction on $p$. The induction base holds by Fact 4. The inductive step follows from the third assertion in Fact 9: $F\leq^{p+1}G$ iff $r(F)\leq^{p}r(G)$ iff  $\bff^n_{m-1}(r(F))\leq^{n+p}\bff^n_{m-1}(r(G))$ iff
$\bff^n_m(F)\leq^{n+p+1}\bff^n_m(G)$.
 \qed

Define  a function $\bff^n:\mathcal{F}_\omega(\bar{k})\to A_k$ by $\bff^n(F)=\bff^n_m(F)$ where $m<\omega$ is the unique number with $F\in\mathcal{F}_m(\bar{k})$, and abbreviate  $\bff^n$ to $\bff$. These function play a prominent role in the next subsection.

\subsection{Main results}\label{main3}

Here we apply the results of the previous subsection to deduce basic facts about the EWH of $k$-partitions in $\mathbb{N}$.

\begin{Theorem}\label{FHalg}
 For any $T\in\mathcal{T}_\omega(\bar{k})$ we have  $\Sigma(\mathbb{N},T)=\{B\in k^\omega\mid B\leq\bff(T)\}$, i.e. $\bff(T)$ is complete in $\Sigma(\mathbb{N},T)$ w.r.t. the reducibility of numberings $\leq$.
\end{Theorem}

{\em Proof.} It suffices to show that for any  $n\geq0$ and any normal  $T\in\mathcal{T}_m(\bar{k})$ we have $\mathcal{L}^{n}(\mathbb{N},T)=\{B\mid B\leq\bff^n_m(T)\}$. Since the base $\mathcal{L}=\{\Sigma^0_{1+n}(\mathbb{N})\}$ is reducible, by Proposition \ref{fineprop}(2) it suffices to check the equality with red-$\mathcal{L}^n(\mathbb{N},T)$ in place of  $\mathcal{L}^n(\mathbb{N},T)$. By induction on $m$ we show that the equality holds for all $n$ and $T$. For $m=0$ the equality holds because both of its parts coincide with $\{\bi\}$ where $T=i<k$. It remains to consider the case $(T,t)\in\mathcal{T}_{m+1}(\bar{k})$, assuming the assertion to hold for $m$.

If $T$ is  singleton  then $\bff^n_{m+1}=\bff^{n+1}_m$ and, by the definition of FH in Section \ref{wkpart}, $\mathcal{L}^{n}(\mathbb{N},T)=\mathcal{L}^{n+1}(\mathbb{N},t(\varepsilon))$ which implies the equality. If $T$ is non-singleton and $p>0$ satisfies $\omega\cap T=\{i\mid i<p\}$ then $\bff^n_{m+1}(T)=\bff^{n+1}_m(t(\varepsilon))\cdot^n(\bigoplus_{i<p}\bff^n_{m+1}(T(i)))$.  For abbreviation, we denote $A=\bff^n_{m+1}(T)$, $\nu=\bff^{n+1}_m(t(\varepsilon))$, and $\mu=\bigoplus_{i<p}\mu_i$ where $\mu_i=\bff^n_{m+1}(T(i))$; thus, we have to show $\mathcal{L}^{n}(\mathbb{N},T)=\{B\mid B\leq A\}$. By induction on the rank of $T$, $\nu$ is complete in $\mathcal{L}^{n+1}(\mathbb{N},t(\varepsilon))$  and $\mu_i$ is complete in $\mathcal{L}^{n}(\mathbb{N},T(i))$ for every $i<p$. By the definition of $\cdot^n$ in Subsection \ref{algebra} and Definition \ref{operG} we have: if $\varphi(r(x))\uparrow$ then $A(x)=\nu(l(x))$ else $A(x)=\mu\varphi(r(x))$, where $\varphi=\tilde{\varkappa}^{\bo^{(n)}}$.

We first check that $A\in\mathcal{L}^{n}(\mathbb{N},T)$. Let $V=\{x\mid y=\varphi(r(x))\downarrow\}$ and $V_i=\{x\mid  rest(y,p)=i\}$ for every $i<p$. Then $V=V_0\cup\cdots\cup V_{p-1}$ and $V_i$ are pairwise disjoint $\Sigma^0_{n+1}$-sets, hence $A=A|_{\overline{V}}\cup A|_{V_0}\cup\cdots\cup A|_{V_{p-1}}$.
For $x\in\overline{V}$ we have $A(x)=\nu(l(x))$ by Proposition \ref{fineprop}(5), hence $A|_{\overline{V}}=(\nu\circ l)|_{\overline{V}}$.
Since $\nu\in\mathcal{L}^{n+1}(\mathbb{N},t(\varepsilon))$, we have $\nu\circ l\in\mathcal{L}^{n+1}(\mathbb{N},t(\varepsilon))$ by Proposition \ref{fineprop}(5), hence $A|_{\overline{V}}\in\mathcal{L}^{n+1}(\overline{V},t(\varepsilon))$ by Proposition \ref{fineprop}(3). 
Since any  $\mu_i$, $i<p$, is a $\bo^{(n)}$-complete numbering and $\varphi\circ r$ is a $\bo^{(n)}$-c.p. function, $\varphi\circ r$ has a computable $\mu_i$-totalizer $g_i$. Thus, $A(x)=\mu_i(g_i(x))$ for each  $x\in V_i$. Arguing as above and taking into account that $\mu_i\in\mathcal{L}^{n}(\mathbb{N},T(i))$, we deduce that  $A|_{V_i}\in\mathcal{L}^{n}(V_i,T(i))$. By Proposition \ref{fineprop}(4), $A\in\mathcal{L}^{n}(\omega,T)$. By Proposition \ref{fineprop}(5),  $\mathcal{L}^{n}(\mathbb{N},T)\supseteq\{B\mid B\leq A\}$. 

It remains to check that any $B\in$red-$\mathcal{L}^{n}(\omega,T)$ is reducible to $A$. Let $B$ be determined by a reduced $T$-family $(\{U_{\tau_0}\}, \{U_{\tau_0\tau_1}\},\ldots)$ in $\mathcal{L}^{n}($. By Proposition \ref{fineprop}(3), $B|_{\overline{U}}\in\mathcal{L}^{n+1}(\overline{U},t(\varepsilon))$ and $B|_{U_i}\in\mathcal{L}^{n}(U_i,T(i))$ for each $i<p$, where $U=U_0\cup\cdots\cup U_{p-1}$. By induction, $B|_{\overline{U}}\leq\nu$ and $B|_{U_i}\leq\mu_i$  for each $i<p$; let $g$ and $g_i$ be  corresponding computable reductions. Let $h(x)=\langle g(x),v(x)\rangle$ where $v$ is a computable function such that $\varphi(v(x))\uparrow$ for $x\not\in U$ and $\varphi(v(x))=p\cdot g_i(x)+i$ for $x\in U_i$, $i<p$. Then $h$ is a computable reduction of $B$ to $A$.   
 \qed
 
Theorem \ref{FHalg}, Fact 4, and Proposition \ref{ncolprop}(3) immediately imply the following result.

\begin{Theorem}\label{main1}
The EWH $\{\Sigma(\mathbb{N},T)\}_{T\in\mathcal{T}_\omega(\bar{k})}$ strongly does not collapse.
\end{Theorem}

Since $\mathbb{N}$ is discrete, the WH $\{\mathbf{\Sigma}(\mathbb{N},T)\}$  collapses to very low levels (it has finitely many distinct levels), the formulation of  Theorem \ref{main1} cannot be improved to the strong uniform version. 

Fact 11 and Proposition \ref{iterh1} imply the following basic result:

\begin{Theorem}\label{alg}
\begin{enumerate}
\item  The function $\bff^n$ induces an isomorphism between the quotient-structure of the structure  $(\mathcal{F}_\omega(\bar{k});\sqcup,\cdot^0,\cdot^{1},\ldots,\leq^0,\leq^1,\ldots)$ under $\equiv_h$ and the quotient-structure of  $(A_k;\oplus,\cdot^n,\cdot^{n+1},\ldots,\leq^n,\leq^{n+1},\ldots)$ under $\equiv^n$. 
\item  The function $\bff=\bff^0$ induces an isomorphism between the quotient-structure  of    $(\mathcal{F}_\omega(\bar{k});\sqcup,\cdot^0,\cdot^{1},\ldots,\leq^0,\leq^1,\ldots)$ under $\equiv_h$ and the quotient-structure of  $(A_k;\oplus,\cdot^0,\cdot^{1},\ldots,\leq^0,\leq^{1},\ldots)$ under $\equiv$. 
\item  The quotient-semilattices of $(A_k;\oplus,\leq)$ and $(\mathcal{F}_\omega(\bar{k});\sqcup,\leq_h)$ are isomorphic.
 \end{enumerate}
\end{Theorem}

The next corollary  would be hard to prove without the established isomorphisms.

\begin{Theorem}\label{compres}
\begin{enumerate}
\item The quotient-structure of  $(A_k;\oplus,\cdot^0,\cdot^{1},\ldots,I,\leq^0,\leq^{1},\ldots)$ under $\equiv$ is computably presentable, where $I$ is the unary relation which is true precisely on  the join-irreducible elements of $(A_k;\oplus,\leq^0)$.
\item All the quotient-semilattice of $(A_k;\oplus,\leq^n)$, $n<\omega$, are isomorphic to each other, and are computably presentable.
 \end{enumerate}
\end{Theorem}

{\em Proof.} (1) For the signature $\{\oplus,\cdot^0,\cdot^{1},\ldots,\leq^0,\leq^{1},\ldots\}$ without $I$, the assertion follows from the previous theorem and the easy fact that the structure  $(\mathcal{F}_\omega(\bar{k});\sqcup,\cdot^0,\cdot^{1},\ldots,\leq^0,\leq^1,\ldots)$ is computably presentable. For the whole signature, it suffices to check that the corresponding relation $I$ of the iterated $k$-forests is computable. This follows from the remarks about minimal forests after Proposition \ref{minq} and at the end of Subsection \ref{trees}.

(2) The isomorphism follows from Fact 11 because all semilattices are isomorphic to the quotient-semilattice of $(\mathcal{F}_\omega(\bar{k});\sqcup,\leq_h)$. The computable presentability follows from (1).
 \qed

The structure $(\{\Sigma(\mathbb{N},T)\mid T\in\mathcal{T}_\omega(\bar{k})\};\subseteq)$ (that is isomorphic to the quotient-poset of $(\mathcal{T}_\omega(\bar{k});\leq_h)$) is rather complicated for $k\geq3$ and very easy (isomorphic to $\bar{2}\cdot\varepsilon_0$) for  $k=2$. It makes sense to look for a natural substructure of $(\{\Sigma(\mathbb{N},T)\mid T\in\mathcal{T}_\omega(\bar{k})\};\subseteq)$ similar to the structure of levels of the FH of sets. In fact, this substructure was identified already in \cite{s83}. We briefly recall it below. 

In the terminology of the present paper, it looks as $(\{\Sigma(\mathbb{N},T_{\alpha,i})\mid i<k,\alpha<\varepsilon_0\};\subseteq)$ where $T_{\alpha,i}=T_{\alpha,i}^0$ and $T^n_{\alpha,i}\in\mathcal{T}_\omega(\bar{k})$ ($n<\omega,i<k,\alpha<\varepsilon_0$) are defined by induction on $\alpha$ as follows: $T^n_{0,i}=i$;  $T^n_{\omega^\gamma,i}=T^{n+1}_{\gamma,i}$ for
$\gamma >0;$
  $T^n_{\beta+1,i}=i\cdot\bigsqcup_{j<k}T^n_{\beta,j}$ for all $\beta<\varepsilon_0$, and
 $T^n_{\beta+\omega^\gamma,i}=T^n_{\omega^\gamma,i}\cdot(\bigsqcup_{j<k}T^n_{\beta,j})$  for $\gamma
>0$ and $\beta$ of the form $\beta =\omega^\gamma\cdot\beta_1 >0$.

\section{EWH in some other spaces}\label{ncol1}

In this section we illustrate the method of Proposition \ref{ncolpres} by proving the non-collapse of EWH and WH in some other concrete spaces. 

We start with observing that Theorem \ref{main1} implies the non-collapse of EWH in some spaces, e.g.:

\begin{Corollary}\label{dunion}
Let $X=X_0\sqcup X_1\sqcup\cdots$ be the disjoint union of a uniform sequence $\{X_n\}$ of nonempty effective cb$_0$-spaces. Then the EWH $\{\Sigma(X,T)\}$ strongly does not collapse.
\end{Corollary}

{\em Proof.} For $x\in X$, let $g(x)$ be the unique number $n$ with $x\in X_n$. Then $g:X\to\mathbb{N}$ is a computable effectively open surjection. By Theorem \ref{main1} and Proposition \ref{ncolpres1}, $\{\Sigma(X,T)\}$ strongly does not collapse.
 \qed

In particular,  Corollary \ref{dunion} applies to $\calN\simeq\calN\sqcup\calN\sqcup\cdots$ which shows that $\{\Sigma(\calN,T)\}$ strongly does not collapse. 
The properties of witnesses for strong non-collapse obtained in this way have very low topological complexity: e.g., the witnesses for the Baire space are clopen (a $k$-partition $A$ is clopen if $A^{-1}(i)$ is clopen for each $i<k$). Therefore, these witnesses cannot be used to show the strong uniform non-collapse property of the EWH in $\mathcal{N}$. We conclude the paper by showing the strong uniform non-collapse of the EWHs in some spaces related to $\mathcal{N}$. These results follow easily from the results in \cite{km19} by observing that their effective versions hold.

First we consider the domain $\bdom$. For this space, the result is obtained by some observations and additions to the proofs in \cite{km19}, so let us recall some information from that paper. 
Let $\home=\omega\cup\{\tt p\}$ be obtained from $\omega$ by adjoining a new element $\rm p$; we endow $\home$ with the discrete topology. For $x\in\conc$, let $\delta(x)\in\bdom$ be obtained from  $x$ by deleting all entries of $\tt p$. A function $f:\conc\to\conc$ (resp.  $A:\conc\to\bar{k}$) is {\em conciliating} if $\delta\circ f=f^*\circ\delta$ (resp. $A=A^*\circ\delta$) for some (unique) $f^*:\bdom\to\bdom$ (resp.  $A^*:\bdom\to\bar{k}$). The function $f$ is {\em initializable} if, for every $\tau\in\home^{<\omega}$, there is a continuous function $h_\tau:\conc\to\conc$ such that $\delta(f(x))=\delta(f(\tau h_\tau(x)))$ for all $x\in\conc$. In Proposition 2.15 from \cite{km19}, an initializable $\bfSig^0_2$-measurable conciliating function $\mathcal{U}:\conc\to\conc$ was constructed which is universal in the sense that for every $\bfSig^0_2$-measurable conciliating function $\mathcal{V}:\conc\to\conc$ there is a continuous function $h:\conc\to\conc$ such that $\delta\circ\mathcal{V}=\delta\circ\mathcal{U}\circ h$.

\begin{Theorem}\label{main}
\begin{enumerate}
 \item The WH $\{\mathbf{\Sigma}(\bdom,T)\}_{T\in\mathcal{T}_\omega(\bar{k})}$ strongly does not collapse. Similarly for the infinitary WH.
 \item The EWH $\{\Sigma(\bdom,T)\}_{T\in\mathcal{T}_\omega(\bar{k})}$ strongly uniformly does not collapse. 
 \end{enumerate}
\end{Theorem}

{\em Proof.} (1) For notation simplicity, we only consider the finitary case, in the infinitary case the argument is the same.   The Definition 3.1.4 in \cite{km19},   using the  induction on trees and the universal function $\mathcal{U}$,  associates  with any tree $T$ a conciliatory $\Omega_T:\conc\to\bar{k}$. By the results in Section 3.3 of \cite{km19}, $\Omega_T$ is in $\mathbf{\Sigma}(\conc,T)\setminus\bigcup\{\mathbf{\Sigma}(\conc,V) \mid V\in\mathcal{T}_\omega(\bar{k}),\;T\not\leq_hV\}$. As the function $\delta$ is a continuous open surjection and $\Omega_T=\Omega_T^*\circ\delta$, we get  $\Omega_T^*\in\mathbf{\Sigma}(\bdom,T)\setminus\bigcup\{\mathbf{\Sigma}(\bdom,V) \mid V\in\mathcal{T}_\omega(\bar{k}),\;T\not\leq_hV\}$ by Proposition \ref{pres}, completing the proof. 

(2) Inspecting the proof of Proposition 2.15 (resp. Lemma 2.11) in \cite{km19} shows that  $\mathcal{U}$ and $\mathcal{U}^*$ are in fact $\Sigma^0_2$-measurable. Inspecting the proof of Lemma 2.15 in \cite{km19} shows that $\Omega_T$ is in $\Sigma(\conc,T)$. Clearly, $\delta$ is a computable effectively open surjection. Thus,  $\Omega_T^*$ is in $\Sigma(\bdom,T)\setminus\bigcup\{\mathbf{\Sigma}(\bdom,V) \mid V\in\mathcal{T}_\omega(\bar{k}),\;T\not\leq_hV\}$ by Proposition \ref{pres}, completing the proof. 
 \qed

This theorem  and Proposition \ref{ncolpres1} imply some new information on the EWH in Baire and Cantor spaces:

\begin{Theorem}\label{main2}
The EWHs $\{\Sigma(\calN,T)\}$ and $\{\Sigma(\mathcal{C},T)\}_{T\in\mathcal{T}_\omega(\bar{k})}$ strongly uniformly do not collapse.
\end{Theorem}

{\em Proof.} As $\delta$ is a computable effectively open sujection,  $\{\Sigma(\conc,T)\}$  strong uniformly does not collapse by Theorem \ref{main} and Proposition \ref{ncolpres}. As $\conc$ is effectively homeomorpic to $\calN$,  the first assertion follows.

For the second assertion, consider the Cantor domain $\cdom$, $n\geq2$, in place of $\bdom$.  A slight modification of the notions from the beginning of this section apply to $\cdom$. Also, a slight modification of the proof of  Theorem \ref{main} shows that it remains true for $\cdom$. As in the previous paragraph, $\{\Sigma(\hat{n}^\omega,T)\}$  strong uniformly does not collapse. Since $\hat{n}^\omega$ is effectively homeomorpic to $(n+1)^\omega$, and thus to $\mathcal{C}$, this implies the second assertion. Note that the spaces $\cdom$ for distinct $n$ are not homeomorphic.
 \qed 

The results and methods above does not directly apply to many important spaces, e.g.  to the Scott domain $P\omega$ and  to the intervals of $\mathbb{R}$. It seems that proving the non-collapse of the (effective) WH in these and many other spaces heavily depends on the topology of a concrete space. Nevertheless, we hope that the methods and results of this paper suggest how one could attack problems of this kind.



%

\end{document}